\documentclass[11pt]{amsart}

\usepackage{geometry}
\usepackage{todonotes}
\usepackage{amsfonts, amssymb, amscd}
\numberwithin{equation}{section}

\usepackage[symbol]{footmisc}

\usepackage{bm}
\usepackage{verbatim}
\usepackage{mathrsfs}
\usepackage{graphicx}
\usepackage{tikz-cd}
\usepackage{subcaption}
\usepackage{listings}
\usepackage{subfiles}
\usepackage[toc,page]{appendix}
\usepackage{mathtools}
\usepackage{comment}
\usepackage{enumerate}
\usepackage{enumitem}
\usepackage[all]{xy}

\usepackage{graphicx}
\graphicspath{{images/}}

\usepackage{appendix}
\usepackage{hyperref}
\hypersetup{
    colorlinks=true,
    citecolor=red,
    linkcolor=blue,
    filecolor=magenta,      
    urlcolor=red,
}
\newcommand{\Mm}{{\bf{M}}}

\newcommand{\Cc}{\mathbb{C}}

\newcommand{\Qq}{\mathbb{Q}}

\newcommand{\Rr}{\mathbb{R}}

\newcommand{\Span}{\operatorname{Span}}

\newcommand{\Center}{\operatorname{center}}

\newcommand{\rk}{\operatorname{rank}}
\newcommand{\red}{\operatorname{red}}

\newcommand{\tang}{{\operatorname{tang}}}

\newcommand{\lct}{\operatorname{lct}}

\newcommand{\fol}{\operatorname{fol}}

\newcommand{\Supp}{\operatorname{Supp}}

\newcommand{\mult}{\operatorname{mult}}

\newcommand{\Ff}{\mathcal{F}}

\newcommand{\Ii}{\Gamma}

\newcommand{\Ee}{\mathcal{E}}

\newcommand{\Sing}{\mathrm{Sing}}

\newtheorem{thm}{Theorem}[section]

\newtheorem{cor}[thm]{Corollary}
\newtheorem{lem}[thm]{Lemma}
\newtheorem{prop}[thm]{Proposition}

\newtheorem{claim}[thm]{Claim}

\theoremstyle{definition}
\newtheorem{defn}[thm]{Definition}

\theoremstyle{definition}
\newtheorem{rem}[thm]{Remark}

\theoremstyle{definition}

\begin{document}

\title{Uniform rational polytopes of foliated threefolds and the global ACC}
\author{Jihao Liu, Fanjun Meng, and Lingyao Xie}

\date{\today}

\begin{abstract}
In this paper, we show the existence of uniform rational lc polytopes for foliations with functional boundaries in dimension $\leq 3$. As an application, we prove the global ACC for foliated threefolds with arbitrary DCC coefficients. We also provide applications on the accumulation points of lc thresholds of foliations in dimension $\leq 3$.
\end{abstract}

\address{Department of Mathematics, Northwestern University, 2033 Sheridan Road, Evanston, IL 60208, USA}
\email{jliu@northwestern.edu}

\address{Department of Mathematics, Johns Hopkins University, 3400 N. Charles Street, Baltimore, MD 21218, USA}
\email{fmeng3@jhu.edu}

\address{Department of Mathematics, The University of Utah, Salt Lake City, UT 84112, USA}
\email{lingyao@math.utah.edu}

\maketitle

\tableofcontents

\section{Introduction}\label{sec:Introduction}

We work over the field of complex numbers $\mathbb{C}$. 

Foliations are important and interesting objects in geometry. Particularly in the minimal model program (MMP), they play a critical role in Miyaoka's proof of some important cases of the abundance conjecture in dimension three \cite{Miy87} (see also \cite[Chapter 9]{Kol+92}). In recent years, progress has been made in the direction of extending the MMP to the setting of foliations, particularly in low dimension. To be specific, the foundations of the MMP for foliated surfaces (cf. \cite{McQ08,Bru15}), foliated threefolds (cf. \cite{CS20,Spi20,CS21,SS22}), and algebraically integrable foliations (cf. \cite{ACSS21,CHLX23,CS23a,LMX24b} have been laid. Many classical questions from the MMP can be asked in the setting of foliations, such as the ascending chain condition (ACC) conjecture for minimal log discrepancies and the ACC conjecture for lc thresholds \cite{Che22,Che23b}.

In a recent work \cite{LLM23}, the first two authors and Y. Luo established the rational case of the global ACC for foliated threefolds, i.e., given any lc foliated log Calabi-Yau triple $(X,\Ff,B)$ of dimension $3$ whose coefficients belong to a set $\Ii$ of rational numbers satisfying the descending chain condition (DCC), the coefficients of $B$ belong to a finite set depending only on $\Ii$. The goal of this paper is to prove the global ACC in its full generality for foliations in dimension $3$. 

\begin{thm}\label{thm: global acc threefold}
Let $\Ii\subset [0,1]$ be a DCC set. Then there exists a finite set $\Ii_0\subset\Ii$ satisfying the following. Let $(X,\Ff,B)$ be a projective lc foliated triple of dimension $\leq 3$ such that $K_{\Ff}+B\equiv 0$ and $B\in\Ii$ (i.e. the coefficients of $B$ belong to $\Ii$). Then $B\in\Ii_0$.
\end{thm}

Theorem \ref{thm: global acc threefold} is an analogue of \cite[Theorem 1.5]{HMX14} for foliated pairs of dimension $\leq 3$. Y.-A. Chen proved Theorem \ref{thm: global acc threefold} in dimension $2$ (\cite[Theorem 2.5]{Che22}).

\begin{rem}[$\mathbb Q$-coefficients versus $\mathbb R$-coefficients]
In an earlier work \cite[Theorem 1.1]{LLM23}, Theorem \ref{thm: global acc threefold} was established for the case when $\Ii\subset\mathbb{Q}$. Despite its technicality, it is important to consider pairs and foliated triple structures with real coefficients for future applications.

For example, when considering the minimal model program with scaling, we need to consider a sequence of scaling numbers $\lambda_i$. In many scenarios we need to consider the limit of the scaling numbers $\lambda:=\lim_{i\rightarrow+\infty}\lambda_i$. However, it is possible that $\lambda$ is an irrational number when each $\lambda_i$ is rational. Thus, we need to study pairs with real coefficients. 

For example, pairs with real coefficients are necessary to develop the theory of \emph{uniform rational polytopes}, which implies the global ACC. The theory of uniform rational polytopes allows us to reduce plenty of questions related to pairs and triples with real coefficients to the case of rational coefficients. Although seemingly technical, this theory (as well as a weaker version, the theory of rational polytopes) is known to be very useful for usual pairs in many different contexts, such as the minimal model program with scaling (\cite{Sho92,BCHM10}), accumulation points of minimal log discrepancies (\cite{Liu18,HLS19,HL22}), construction of special minimal model programs (\cite{Kol21,HL23}), and the effective Iitaka fibration conjecture \cite{CHL23}. For further details about this theory, we refer the reader to \cite{Nak16,HLS19,HLQ21}.
\end{rem}

Due to the failure of effective birationality for foliations \cite[Paragraph before Theorem 1.4]{SS23}, the proof of the rational case of the global ACC for foliated threefolds heavily relies on the index theorem of surfaces, which only works for the rational coefficient case. Therefore, there are essential difficulties in extending \cite[Theorem 1.1]{LLM23} to the general case when the coefficients of $B$ belong to an arbitrary DCC set $\Ii\subseteq\mathbb{R}$. 

To resolve this issue, we prove the existence of uniform rational lc polytopes for foliations in dimension $\leq 3$. This result is crucial for the proof of our main theorem, Theorem \ref{thm: global acc threefold}, and has potential applications in the study of foliations and their singularities.

\begin{thm}\label{thm: uniform rational polytope foliation intro}
Let $v_1^0,\dots,v_m^0$ be positive numbers and $\bm{v}_0:=(v_1^0,\dots,v_m^0)$. Then there exists an open set $U\ni \bm{v}_0$ of the rational envelope of $\bm{v}_0$ satisfying the following. 

Let $(X,\Ff,B=\sum_{j=1}^mv_j^0B_j)$ be a projective foliated lc triple of dimension $\leq 3$, where $B_j\geq 0$ are distinct Weil divisors. Then  $(X,\Ff,B=\sum_{j=1}^mv_jB_j)$ is lc for any $(v_1,\dots,v_m)\in U$.
\end{thm}
We remark that the projectivity condition in Theorem \ref{thm: uniform rational polytope foliation intro} is expected to be unnecessary. We add this condition for consistency with references, e.g. \cite{CS20}.

When $\dim X=2$, Theorem \ref{thm: uniform rational polytope foliation intro} was proved in \cite{LMX24a}. Theorem \ref{thm: uniform rational polytope foliation intro} is an analogue of the theory of uniform rational lc polytopes for usual pairs with functional boundaries \cite[Theorem 5.6]{HLS19}. 

The proof of Theorem \ref{thm: global acc threefold} is an application of the theory of uniform rational polytopes which allows us to reduce Theorem \ref{thm: global acc threefold} to the case of rational coefficients. This reduction greatly simplifies the arguments and allows us to bypass the difficulties arising from the lack of effective birationality for foliations.

As a direct corollary of Theorem \ref{thm: uniform rational polytope foliation intro}, we obtain the following result on the accumulation points of lc thresholds of foliations. Again, the projectivity condition is added in Corollary \ref{cor: accumulation point foliated threefold lct} only for consistency with references.

\begin{cor}\label{cor: accumulation point foliated threefold lct}
Let $\Ii\subset [0,1]$ be a DCC set such that $\bar\Ii\subset\mathbb Q$. We consider the set of foliated lc thresholds in dimension $3$
$$\lct_{\fol}(3,\Ii):=\{\lct(X,\Ff,B;D)\mid \dim X=3, X\text{ is projective, } B\in\Ii, D\in\mathbb N^+\}.$$
Then the accumulation points of $\lct_{\fol}(3,\Ii)$ are rational numbers.
\end{cor}

\noindent\textit{Idea and sketch of the proof}. We start with the proof of Theorem \ref{thm: global acc threefold}. Since we work in dimension $3$, it is not difficult to reduce Theorem \ref{thm: global acc threefold} to the case when there exists a contraction $\pi: X\rightarrow Z$ such that the general fibers of $\pi$ are tangent to $\Ff$ and $\dim Z>0$, i.e. the setting of Proposition \ref{prop: global acc with contraction}. In this case, the usual global ACC implies that the coefficients of the horizontal$/Z$ part of $B$ belong to a finite set, so we only need to worry about the vertical$/Z$ part of $B$. We let $B_Z$ be the discriminant part of the canonical bundle formula of $\pi: (X,\Ff,B)\rightarrow Z$. Then the coefficients of $B_Z$ belong to a DCC set by the ACC for lc thresholds of foliations.

At this point, we prove by using contradiction: if the coefficients of $B$ do not belong to a finite set, then we may construct a boundary $\bar B\geq B$ sufficiently close to $B$ whose coefficients belong to a finite set, and $K_{\Ff}+\bar B\sim_{\mathbb R,Z}0$. By the ACC for lc thresholds of foliations, the discriminant part $\bar B_Z$ induced by the canonical bundle formula of $\pi: (X,\Ff,\bar B)\rightarrow Z$ has coefficients which are larger and sufficiently close to the ones of $B_Z$, and its coefficients also belong to a DCC set. Therefore, the coefficients of $B_Z$ belong to a DCC but not finite set. For the rational coefficient case, we get a contradiction to the global ACC of foliated log Calabi-Yau triples polarized with a semi-ample divisor \cite[Lemma 7.2]{LLM23}, but this can no longer work for the real coefficient case. 

To deal with this issue, the key idea is the following: we would like to find positive real numbers $a_1,\dots,a_k$ depending only on the coefficient set $\Ii$, so that we have a decomposition $K_{\Ff}+\bar B=\sum a_i(K_{\Ff}+\bar B_i)$, where $\sum a_i=1$, $(X,\Ff,\bar B_i)$ is lc, $\bar B_i$ is a $\Qq$-divisor, and $K_{\Ff}+\bar B_i\sim_{\mathbb Q,Z}0$ for each $i$. With such decomposition established, we may let $\bar B_{Z_i}$ be the discriminant part induced by the canonical bundle formula of $\pi: (X,\Ff,\bar B_i)\rightarrow Z$ for each $i$. By using similar arguments as the rational coefficient case, we can compare the coefficients of $\sum a_i\bar B_{Z_i}$ and $B_Z$ and get a contradiction, which leads to a proof of Theorem \ref{thm: global acc threefold}. The remaining difficulty is to find such a decomposition $K_{\Ff}+\bar B=\sum a_i(K_{\Ff}+\bar B_i)$, which is nothing but Theorem \ref{thm: uniform rational polytope foliation intro}. 

Next we turn to the proof of Theorem \ref{thm: uniform rational polytope foliation intro}. Due to technicality, in the following, we provide the reader with a sketch of the proof of Theorem \ref{thm: uniform rational polytope foliation intro} for the following special case: $m=1$, $\bm{v}_0=\left(\frac{\sqrt{2}}{2}\right)$, and $S:=B_1$ is a prime divisor. The general case follows from similar lines of proof. 

To prove this special case of Theorem \ref{thm: uniform rational polytope foliation intro}, we only need to show that there exists a rational number $a>\frac{\sqrt{2}}{2}$ (which does not depend on $X$) such that $(X,\Ff,aS)$ is lc. By passing to a dlt model of $(X,\Ff,\frac{\sqrt{2}}{2}S)$ we may assume that $(X,\Ff,\frac{\sqrt{2}}{2}S)$ is $\Qq$-factorial dlt. (Here divisors with coefficient $1$ may appear but will not influence the proof. For simplicity, let us ignore them.) Suppose Theorem \ref{thm: uniform rational polytope foliation intro} does not hold in this case, then there exists a sequence $(X_i,\Ff_i,a_iS_i)$ such that $a_i$ is strictly decreasing, $\lim_{i\rightarrow+\infty}a_i=\frac{\sqrt{2}}{2}$, and $a_i$ is the lc threshold of $S_i$ with respect to $\Ff_i$. We suppose that $E_i$ is a prime divisor which achieves the lc threshold. Since $(X_i,\Ff_i,\frac{\sqrt{2}}{2}S_i)$ is $\Qq$-factorial dlt, it can be reduced to the case when $E_i$ is $\Ff_i$-invariant. Moreover, we may assume that $\Center_{X_i}E_i$ is a closed point, as other cases can be reduced to the lower dimensional cases which are proved in \cite{LMX24a}. 

Now there is a key observation: suppose that there exists an extraction $g_i: Y_i\rightarrow X_i$ of $E_i$. Let $\Ff_{Y_i}:=g_i^{-1}\Ff_i$ and $S_{Y_i}:=(g_i^{-1})_*S_i$, then $K_{\Ff_{Y_i}}+a_iS_{Y_i}\sim_{\mathbb R,X_i}0$ but $K_{\Ff_{Y_i}}+\frac{\sqrt{2}}{2}S_{Y_i}$ is anti-ample$/X_i$. Therefore, for the normalization $V_i$ of any non-trivial lc center of $(Y_i,\Ff_{Y_i},a_iS_{Y_i})$ in $E_i$ (e.g. normalization of $E_i$), $(K_{\Ff_{Y_i}}+a_iS_{Y_i})|_{V_i}\equiv 0$ but $(K_{\Ff_{Y_i}}+\frac{\sqrt{2}}{2}S_{Y_i})|_{V_i}\not\equiv 0$. Therefore, if we can achieve precise adjunction formulas of foliated threefolds to $\Ff_{Y_i}$-invariant lc centers, the coefficients of $(K_{\Ff_{Y_i}}+a_iS_{Y_i})|_{V_i}$ and $(K_{\Ff_{Y_i}}+\frac{\sqrt{2}}{2}S_{Y_i})|_{V_i}$ will be controlled. Now we can use lower dimensional results to get a contradiction in this case. With this observation in mind, there are two things we need to do:
\begin{enumerate}
    \item Establish precise adjunction formulas to invariant lc centers. The whole Section \ref{sec: adjunction} is dedicated to it. The proofs generally follow from the same lines of the proofs of the adjunction formulas in \cite{CS20,CS21} but we provide more details. It is important to note that for the rank two case, when the minimal lc center is not a divisor, we need to consider adjunction to the normalization of lc centers of dimension $1$ as well (see Theorem \ref{thm: adjunction formula corank 1 dim 3 intro}(6)).
    \item Show the existence of the extraction $g_i: Y_i\rightarrow X_i$ of $E_i$. This is our key lemma, see Lemma \ref{lem: find nontrivial divisor on dlt model}, and Subsections \ref{subsec: lct}, \ref{subsec: special dlt} are dedicated to it. More precisely, we shall prove the existence of such $g_i: Y_i\rightarrow X_i$ after a suitable substitution of $X_i$ and $E_i$.
    
    The proof of Lemma \ref{lem: find nontrivial divisor on dlt model} is as follows: first we take a dlt modification $W_i\rightarrow X_i$ of $(X_i,\Ff_i,a_iS_i)$ which extracts $E_i$. It may extract some divisors that are lc centers of $(X_i,\Ff_i,\frac{\sqrt{2}}{2}S_i)$. By replacing $X_i$ with a higher model, we may assume that any divisor extracted by $W_i\rightarrow X_i$ is not an lc place of $(X_i,\Ff_i,\frac{\sqrt{2}}{2}S_i)$. Now $K_{\Ff_{W_i}}+\frac{\sqrt{2}}{2}S_{W_i}\sim_{\mathbb R,X_i}F_i$ where $F_i\geq 0$ is supported on the exceptional divisor of $W_i\rightarrow X_i$. Now we may run a $(K_{\Ff_{W_i}}+\frac{\sqrt{2}}{2}S_{W_i})$-MMP$/X_i$ with scaling of an ample divisor. By the general negativity lemma \cite[Lemma 3.3]{Bir12} this MMP contracts $F_i$. Since $X_i$ is $\Qq$-factorial, it must terminate with $X_i$, and the last step of the MMP must be a divisorial contraction. We may let  $g_i: Y_i\rightarrow X_i$  be the last step of this MMP. By our construction, $g_i$ extracts a prime divisor $E_i'$ which achieves the lc threshold of $S_i$ with respect to $\Ff_i$ (although it is possible that $E_i'\not=E_i$). This suffices our requirements.

    It is worth to mention that our key Lemma, Lemma \ref{lem: find nontrivial divisor on dlt model}, is proved to be important in the study of log canonicity for foliations, although we cannot find any analogue of it in the study of log canonicity for usual pairs. In later works, analogues of Lemma \ref{lem: find nontrivial divisor on dlt model} have played crucial roles in the proofs of the ACC for lc thresholds and the existence of uniform rational polytopes for algebraically integrable foliations (\cite[Lemma 5.1]{DLM23}) and algebraically integrable generalized foliated quadruples (\cite[Lemma 10.2.1]{CHLX23}).
\end{enumerate}
Finally, with both (1) and (2) settled, we can get pair (resp. foliated pair) structures of $(K_{\Ff_{Y_i}}+a_iS_{Y_i})|_{V_i}$ and $(K_{\Ff_{Y_i}}+\frac{\sqrt{2}}{2}S_{Y_i})|_{V_i}$ with coefficients controlled and good numerical properties when $\rk\Ff=2$ (resp. $\rk\Ff=1$). Now we may use \cite[Theorem 3.6]{Che23a} and \cite[Theorem 3.8, Corollary 3.9]{Nak16} to study the coefficients of $(K_{\Ff_{Y_i}}+a_iS_{Y_i})|_{V_i}$ and $(K_{\Ff_{Y_i}}+\frac{\sqrt{2}}{2}S_{Y_i})|_{V_i}$ and conclude our proof. 

To be specific, when $\rk\Ff=2$, \cite[Theorem 3.6]{Che23a} and \cite[Theorem 3.8, Corollary 3.9]{Nak16} can be applied directly to our setting to check the coefficients of $(K_{\Ff_{Y_i}}+a_iS_{Y_i})|_{V_i}$ and $(K_{\Ff_{Y_i}}+\frac{\sqrt{2}}{2}S_{Y_i})|_{V_i}$ and conclude the proof of Theorem \ref{thm: uniform rational polytope foliation intro} in this case (see Subsection \ref{subsec: 3fold rk2}). When $\rk\Ff=1$, it is easy to check that the restricted foliation $\Ff_{E_i^\nu}$ of $\Ff_i$ on the normalization $E_i^\nu$ of $E_i$ has a non-pseudo-effective foliated canonical divisor, hence it is algebraically integrable. Therefore, we can consider the intersection numbers of $(K_{\Ff_{Y_i}}+a_iS_{Y_i})|_{E^\nu_i}$ and $(K_{\Ff_{Y_i}}+\frac{\sqrt{2}}{2}S_{Y_i})|_{E^\nu_i}$ with a general member of the family of leaves of $\Ff_{E_i^\nu}$ and apply \cite[Theorem 3.6]{Che23a} and \cite[Theorem 3.8, Corollary 3.9]{Nak16} again to conclude the proof of Theorem \ref{thm: uniform rational polytope foliation intro} in this case (see Subsection \ref{subsec: 3fold rk1}).

\medskip

\noindent\textit{Structure of this paper.} Section \ref{sec: Preliminaries} introduces some preliminary results for foliations. Section \ref{sec: adjunction} proves the precise adjunction formulas for foliations in dimension $\leq 3$ to invariant divisors (Theorems \ref{thm: adjunction surface intro}, \ref{thm: adjunction threefold rank 1 intro}, and \ref{thm: adjunction formula corank 1 dim 3 intro}). Section \ref{sec: Uniform rational polytopes} establishes the theory of uniform rational polytopes for foliations with functional boundaries in dimension $\leq 3$, and proves Theorem \ref{thm: uniform rational polytope foliation intro}. Section \ref{sec: proof of the main theorems} provides the proofs of Theorem \ref{thm: global acc threefold} and Corollary \ref{cor: accumulation point foliated threefold lct}.

\medskip

\noindent\textbf{Acknowledgements}. We thank Paolo Cascini, Yen-An Chen, Omprokash Das, Christopher D. Hacon, Jingjun Han, and Yuchen Liu for helpful discussions. We would like to acknowledge the assistance of ChatGPT in polishing the wording.  The third author is partially supported by NSF research grants no: DMS-1801851, DMS-1952522 and by a grant from the Simons Foundation; Award Number: 256202.

\section{Preliminaries}\label{sec: Preliminaries}

We will work over the field of complex numbers $\Cc$. Throughout the paper, we will mainly work with normal quasi-projective varieties to ensure consistency with the references. However, most results should also hold for normal varieties that are not necessarily quasi-projective. Similarly, most results in our paper should hold for any algebraically closed field of characteristic zero. We will adopt the standard notations and definitions in \cite{Sho92,KM98, BCHM10} and use them freely. For foliations and foliated (sub-)triples, we will follow the notations and definitions in \cite{LLM23}, which are mostly consistent with those in \cite{CS20,Spi20,ACSS21,CS21}. For generalized pairs and generalized foliated quadruples, we will follow the notations and definitions in \cite{HL23,LLM23}.

\subsection{Sets}

\begin{defn}\label{defn: DCC and ACC}
Let $\Ii\subset\Rr$ be a set. We say that $\Ii$ satisfies the \emph{descending chain condition} (DCC) if any decreasing sequence in $\Ii$ stabilizes, and $\Ii$ satisfies the \emph{ascending chain condition} (ACC) if any increasing sequence in $\Ii$ stabilizes. We define
$$\Ii_+:=\left\{0\right\}\cup\left\{\sum_{i=1}^n \gamma_i\biggm| n\in\mathbb N^+, \gamma_1,\dots,\gamma_n\in\Ii\right\}.$$
and define $\bar\Ii$ to be the closure of $\Ii$ in $\mathbb R$.
\end{defn}

\begin{defn}
Let $m$ be a positive integer and $\bm{v}\in\mathbb R^m$. The \emph{rational envelope} of $\bm{v}$ is the minimal rational affine subspace of $\mathbb R^m$ which contains $\bm{v}$. For example, if $m=2$ and $\bm{v}=(\frac{\sqrt{2}}{2},1-\frac{\sqrt{2}}{2})$, then the rational envelope of $\bm{v}$ is $(x_1+x_2=1)\subset\mathbb R^2_{x_1x_2}$.
\end{defn}

\subsection{Foliations}

\begin{defn}[Special divisors on foliations, cf. {\cite[Definition 2.2]{CS21}}]\label{defn: special divisors on foliations}
Let $X$ be a normal variety and $\Ff$ a foliation on $X$. For any prime divisor $C$ on $X$, we define $\epsilon_{\Ff}(C):=1$ if $C$ is not $\Ff$-invariant, and  $\epsilon_{\Ff}(C):=0$ if $C$ is $\Ff$-invariant. If $\Ff$ is clear from the context, then we may use $\epsilon(C)$ instead of $\epsilon_{\Ff}(C)$. For any $\Rr$-divisor $D$ on $X$, we define $$D^{\Ff}:=\sum_{C\mid C\text{ is a component of }D}\epsilon_{\Ff}(C)C.$$
Let $E$ be a prime divisor over $X$ and $f: Y\rightarrow X$ a projective birational morphism such that $E$ is on $Y$. We define $\epsilon_{\Ff}(E):=\epsilon_{f^{-1}\Ff}(E)$. It is clear that $\epsilon_{\Ff}(E)$ is independent of the choice of $f$.
\end{defn}

\begin{defn}\label{defn: log resolution}
Let $X$ be a normal variety, $\Ff$ a foliation on $X$, and $B$ an $\Rr$-divisor on $X$, such that either $\rk\Ff=\dim X-1$ or $\dim X=3$ and $\rk\Ff=1$. A \emph{foliated log resolution} of $(X,\Ff,B)$ is a projective birational morphism $f: Y\rightarrow X$ such that $(Y,\Ff_Y:=f^{-1}\Ff,B_Y:=f^{-1}_*B+E)$ is foliated log smooth (cf. {\cite[Definition 3.1]{CS21}, \cite[Definition 4.2]{LLM23}}), where $E$ is the reduced exceptional divisor of $f$. A \emph{foliated resolution} of $\Ff$ is a foliated log resolution of $(X,\Ff,0)$.
\end{defn}

\subsection{Dlt models}

\begin{defn}[Dlt singularities]\label{def: fdlt}
Let $(X,\Ff,B)$ be a foliated triple. 
\begin{enumerate}
    \item Suppose that $\rk\Ff=\dim X-1$. We say that $(X,\Ff,B)$ and $(\Ff,B)$ are \emph{dlt} if 
\begin{enumerate}
\item $(X,\Ff,B)$ is lc,
\item every component of $B$ is generically transverse to $\Ff$, and
\item there exists a foliated log resolution of $(X,\Ff,B)$ which only extracts divisors $E$ of discrepancy $>-\epsilon(E)$.
\end{enumerate}
\item Suppose that $\dim X=3$ and $\rk\Ff=1$.  We say that $(X,\Ff,B)$ and $(\Ff,B)$ are \emph{dlt} if 
\begin{enumerate}
\item $(X,\Ff,B)$ is lc, and
\item $\Ff$ has simple singularities \cite[Definition 2.32]{CS20}.
\end{enumerate}
\end{enumerate}
We remark that dlt implies non-dicritical (cf. \cite[Theorem 11.3]{CS21} and \cite[Lemma 2.8]{CS20}).
\end{defn}

\begin{defn}[Dlt modification]\label{defn: fdlt modification}
Let $(X,\Ff,B)$ be an lc foliated triple. An \emph{dlt modification} of $(X,\Ff,B)$ is a birational morphism $f: Y\rightarrow X$ satisfying the following. Let $\Ff_Y:=f^{-1}\Ff$, $E$ the reduced exceptional divisor of $f$, and  $B_Y:=f^{-1}_*B+E^{\Ff_Y}$.
\begin{enumerate}
    \item $Y$ is $\mathbb Q$-factorial klt, 
    \item $K_{\Ff_Y}+B_Y=f^*(K_{\Ff}+B)$, and
    \item $(Y,\Ff_Y,B_Y)$ is dlt.
\end{enumerate}
We say that $(Y,\Ff_Y,B_Y)$ is a \emph{dlt model} of $(X,\Ff,B)$ and $(\Ff_Y,B_Y)$ is a \emph{dlt model} of $(\Ff,B)$.
\end{defn}

\subsection{A perturbation formula}

We recall a weaker version of Theorem \ref{thm: uniform rational polytope foliation intro} which will be crucial for the proof of our theorems.

\begin{thm}[{\cite[Theorem 1.7]{LLM23}}]\label{thm: rational polytope foliation intro}
Let $v_1^0,\dots,v_m^0$ be positive integers, $\bm{v}_0:=(v_1^0,\dots,v_m^0)$, and $(X,\Ff,B=\sum_{i=1}^mv_i^0B_i)$ an lc foliated triple of dimension $\leq 3$, where $B_i\geq 0$ are distinct Weil divisors. Then there exists an open set $U\ni \bm{v}_0$ of the rational envelope of $\bm{v}_0$, such that $(X,\Ff,B=\sum_{i=1}^mv_iB_i)$ is lc for any $(v_1,\dots,v_m)\in U$.
\end{thm}
Note that the major difference between Theorem \ref{thm: rational polytope foliation intro} and Theorem \ref{thm: uniform rational polytope foliation intro} is that the set $U$ may depend on $(X,\Ff,B)$ in Theorem \ref{thm: rational polytope foliation intro}.

\section{Precise adjunction formulas to invariant divisors}\label{sec: adjunction}

To prove Theorem \ref{thm: uniform rational polytope foliation intro}, we need several precise adjunction formulas to divisors that are invariant to foliations. Although there are many results on adjunction formulas of foliations in the literature \cite{Bru02,Bru15,CS20,Spi20,CS21,SS22,Che22}, the statements of these results are insufficient for our purposes, as we need an accurate description of the coefficients of the different foliated triples. For instance, when applying the adjunction formula to two distinct foliated triples $(X,\Ff,B=\sum b_jB_j)$ and $(X,\Ff,B'=\sum b_j'B_j)$, we need to understand the relationship between the coefficients of the differents of two distinct triples. Additionally, we require adjunction formulas for triples with $\Rr$-coefficients rather than $\Qq$-coefficients. In this section, we provide precise adjunction formulas for invariant divisors, even if they may be well-known to experts. The key differences between these formulas and the traditional adjunction formulas for foliations are the following:
\begin{itemize}
    \item We have a more accurate control of the coefficients of the foliated different $C_i$. This will be helpful when applying adjunction to two or more triples and considering their behavior.
    \item We deal with $\Rr$-coefficients rather than $\Qq$-coefficients.
    \item We only need to control the singularities of one foliated triple $(X,\Ff,B=\sum b_jB_j)$ to get an adjunction formula for all foliated triples of the form $(X,\Ff,B'=\sum b_j'B_j)$, even if the singularities of the latter triple may not be controlled.
\end{itemize}
We refer the reader to \cite[Lemma 3.18]{CS21}, \cite[Proposition 2.16]{CS20}, \cite[Lemma 3.11]{SS22}, \cite[Proposition 3.10, Theorem 4.6]{Che22} for related references. 

Finally, we remark that a recent paper of Cascini and Spicer \cite{CS23b} has established the adjunction formula of foliations to non-invariant divisors in all dimensions. We expect that the ideas in \cite{CS23b} to provide a precise adjunction formula of foliations to non-invariant divisors in any dimensions.

\subsection{Surface case}

In this subsection we prove Theorem \ref{thm: adjunction surface intro}. We actually do not need the precise adjunction formula for surfaces to prove our main theorems, Theorems \ref{thm: global acc threefold} and \ref{thm: uniform rational polytope foliation intro}, but we include the result for completeness.

\begin{lem}\label{lem: small multiplicities surface}
Let $(X\ni x,\Ff,B)$ be an lc (resp. terminal) foliated germ such that $X\ni x$ is terminal. Let $f: Y\rightarrow X$ be the minimal resolution of $X\ni x$ and $I:=\det(\mathcal{D}(f))$. Then there exists a unique $\Ff$-invariant curve $C$ which passes through $x$, and $(B\cdot C)_x\leq\frac{1}{I}$ (resp. $<\frac{1}{I}$).
\end{lem}
\begin{proof}
The existence and uniqueness of $C$ follows from \cite[Theorems 3.19, 4.1]{LMX24a}. Moreover, by \cite[Theorem 4.1]{LMX24a}, $(X\ni x,B+C)$ is lc. Let
$$K_{C^\nu}+B_C:=(K_X+B+C)|_{C^{\nu}},$$
where $C^\nu$ is the normalization of $C$, then $(C^\nu\ni x,B_C)$ is lc. By \cite[Theorem 3.19]{LMX24a}, $I$ is the local Cartier index of $X\ni x$. Thus
$$1\geq\text{(resp. }>\text{)}\mult_xB_C=\frac{I-1}{I}+(B\cdot C)_x$$
and the lemma follows.
\end{proof}

\begin{thm}[Surface case]\label{thm: adjunction surface intro}
Let $(X,\Ff,B=\sum_{j=1}^mb_jB_j)$ be a dlt foliated triple such that $\dim X=2,\rk\Ff=1$, and $B_j$ are the irreducible components of $B$. Assume that
\begin{itemize}
    \item  $C$ is an $\Ff$-invariant curve on $X$,
    \item  $C^\nu$ is the normalization of $C$,
    \item  $P_1,\dots,P_n$ are all closed points on $C$ that are singular points of $X$ or contained in $C\cap\Supp B$, and
    \item $Q_1,\dots,Q_l$ are all closed points on $C$ such that $Z(\Ff,C,Q_i)\in\mathbb N^+$. We refer the reader to \cite[Section 2]{Bru02} for the definition of $Z(\Ff,C,Q_i)$.
\end{itemize}
Then there exist positive integers $w_1,\dots,w_n$ and non-negative integers $\{w_{i,j}\}_{1\leq i\leq n,1\leq j\leq m}$, such that for any real numbers $b_1',\dots,b_m'$, we have the following.
\begin{enumerate}
    \item By identifying $Q_i$ with its inverse image in $C^\nu$ under the normalization $C^\nu\rightarrow C$, we have
    $$\left(K_{\Ff}+\sum_{j=1}^mb_j'B_j\right)\Bigg|_{C^\nu}=K_{C^\nu}+\sum_{i=1}^n\frac{w_i-1+\sum_{j=1}^mw_{i,j}b_j'}{w_i}P_i+\sum_{i=1}^lZ(\Ff,C,Q_i)Q_i.$$
    \item $w_i$ is the local Cartier index (i.e. order of the local fundamental group) of $P_i$ for each $i$.
    \item For any $i$ such that $w_i=1$, $w_{i,j}>0$ for some $j$.
    \item If $(X,\Ff,\sum_{j=1}^mb_j'B_j)$ is lc, then $$\left(C^\nu,\sum_{i=1}^n\frac{w_i-1+\sum_{j=1}^mw_{i,j}b_j'}{w_i}P_i+\sum_{i=1}^lQ_i\right)$$
    is lc, i.e. $\sum_{j=1}^mw_{i,j}b_j'\leq 1$ for any $i$.
\end{enumerate}
\end{thm}

\begin{proof}
Since $(X,\Ff,B)$ is dlt, $X$ is klt (cf. \cite[Corollary 3.20]{LMX24a}), so $X$ is $\Qq$-factorial.

First we prove (1)(2). By definition, for any closed point $Q\not\in\Supp(B)\cup\Sing(X)$, the vanishing order of $K_{\Ff}|_{C^\nu}$ at $Q$ is equal to $Z(\Ff,C,Q)$. Thus if $Q$ is a non-singular point of $C$, then we may identify $Q$ with its inverse image in $C^\nu$. If $Q$ is a singular point of $C$, then since $(X,\Ff,B)$ is dlt, $Q$ is a nodal singularity of $C$ and $Z(\Ff,C,Q)=0$. 

Therefore, to prove (1)(2), we only need to show that for any $1\leq i\leq n$, $C$ is non-singular at $P_i$ and the vanishing order of $(K_{\Ff}+\sum_{j=1}^mb_j'B_j)|_{C^\nu}$ at $P_i$ is equal to $\frac{w_i-1+\sum_{j=1}^nw_{i,j}b_j'}{w_i}$, and $w_i$ is equal to the index of $P_i$.

By \cite[Theorem 3.19]{LMX24a}, each $P_i$ is either a non-singular point of $X$ or a cyclic quotient singularity of $X$. By \cite[Theorem 4.1]{LMX24a}, $(X\ni P_i,C)$ is plt. Thus $C$ is non-singular at $P_i$.

We let $\pi_i: X_i\rightarrow\tilde X$ be an index $1$ cover of $X\ni P_i$, $P_i':=\pi_i^{-1}(P_i)$, $C_i:=\pi_i^*C$, and $B_{j,i}:=\pi_i^*B_j$ for any $j$. Then $C_i$ is non-singular at $P_i'$, and $$\mult_{P_i'}\left(\left(\sum_{j=1}^mb_j'B_{j,i}\right)\Biggm|_{C_i}\right)=\sum_{j=1}^mb_j'(B_{j,i}\cdot C_i)_{P_i'}.$$
Moreover, we have $\deg(\pi_i)=w_i$, the index of $P_i$, for each $i$. We let $w_{i,j}:=(B_{j,i}\cdot C_i)_{P_i'}$ for any $i,j$. (1)(2) now follow from the Hurwitz formula (cf. \cite[Proposition 3.7]{Spi20}). 

Since $P_i\in\Sing(X)\cup(C\cap\Supp B)$, if $w_i=1$, then $P_i\in C\cap\Supp B$. Thus $(B_{j,i}\cdot C_i)_{P_i'}\not=0$ for some $j$, and (3) follows.

Since $w_{i,j}:=(B_{j,i}\cdot C_i)_{P_i'}$ for any $i,j$, by Lemma \ref{lem: small multiplicities surface},
$$\frac{1}{w_i}\geq (B\cdot C)_{P_i}=\frac{1}{w_i}(\pi_i^*B\cdot C_i)_{P_i'}=\frac{1}{w_i}
\left(\sum_{j=1}^mb_j'B_{j,i}\cdot C_i\right)_{P_i'}=\sum_{j=1}^m\frac{w_{i,j}b_j'}{w_i},$$
and (4) follows. 
\end{proof}

\subsection{Threefold rank one case}

\begin{thm}[Threefold rank one case]\label{thm: adjunction threefold rank 1 intro}
Let $(X,\Ff,B=\sum_{j=1}^mb_jB_j)$ be a dlt foliated triple such that $\dim X=3,\rk\Ff=1$, and $B_j$ are the irreducible components of $B$. Let $S$ be an $\Ff$-invariant surface in $X$ and $S^\nu$ the normalization of $S$. Suppose that $B_j$ is $\Qq$-Cartier for any $j$. Then there exist a rank $1$ foliation $\Ff_S$ on $S^\nu$, prime divisors $C_1,\dots,C_n$ on $S^\nu$, and non-negative integers $\{w_{i,j}\}_{1\leq i\leq n, 0\leq j\leq m}$, such that for any real numbers $b_1',\dots,b'_m$, we have the following adjunction formula
$$\left(K_{\Ff}+\sum_{j=1}^mb_j'B_j\right)\Bigg|_{S^\nu}=K_{\Ff_S}+\sum_{i=1}^n\frac{w_{i,0}+\sum_{j=1}^mw_{i,j}b_j'}{2}C_i.$$
We remark that there is no control of the singularities of $\left(S^\nu,\Ff_S,\sum_{i=1}^n\frac{w_{i,0}+\sum_{j=1}^mw_{i,j}b_j'}{2}C_i\right)$ even if $(X,\Ff,\sum_{j=1}^mb_j'B_j)$ is dlt.
\end{thm}

\begin{proof}
The proof is parallel to \cite[Propositon 2.16]{CS20}. For the reader's convenience we give a full proof here. 

By \cite[Proposition-Definition 3.6, Remark 3.7]{CS23b}, there exists a naturally defined restricted foliation $\Ff_S$ on $S^\nu$ of $\Ff$ such that  $\rk\Ff_S=1$. We prove the following claim.

\begin{claim}\label{claim: 2kf cartier codimension 2}
For any prime divisor $D$ on $S^\nu$ and any prime divisor $L$ on $X$, $2L$ is Cartier near the generic point of the image of $D$ in $S$.
\end{claim}
\begin{proof}
Since $(X,\Ff,B)$ is dlt, $\Ff$ has simple singularities. Let $\eta_D$ be the generic point of the image of $D$ in $S$. If $\Ff$ has canonical but not terminal singularity near $\eta_D$, then by the definition of simple singularities, $2L$ is Cartier near $\eta_D$. Thus we may assume that $\Ff$ is terminal near $\eta_D$. By \cite[Lemma 2.12]{CS20}, the image of $D$ in $S$ is not $\Ff$-invariant. By \cite[Lemma 2.6]{CS20}, $\eta_D$ is not contained in $\Sing(X)$, so $L$ is Cartier near $\eta_D$ and we are done.
\end{proof}

\noindent\textit{Proof of Theorem \ref{thm: adjunction threefold rank 1 intro} continued}. By Noetherian property, we only need to show that for any prime divisor $D$ on $S^\nu$, there exist non-negative integers $d_0,\dots,d_n$ such that
$$\left(K_{\Ff}+\sum_{j=1}^mb_j'B_j\right)\Biggm|_{S^\nu}=K_{\Ff_S}+\frac{d_0+\sum_{j=1}^md_jb_j'}{2}D$$
for any real numbers $b_1',\dots,b_m'$ near the generic point of $D$. By Claim \ref{claim: 2kf cartier codimension 2}, $2K_{\Ff}$ and $2B_j$ are Cartier near the image of the generic point of $D$ on $S$. Then we may let $d_j:=2\mult_D(B_j|_{S^\nu})$ for any $j\geq 1$. Moreover, the map $\Omega^1_{X}\otimes\Omega^1_{X}\rightarrow\mathcal{O}_X(2K_{\Ff})$ naturally restricts to a map  $\Omega^1_{X}|_S\otimes\Omega^1_{X}|_S\rightarrow\mathcal{O}_S(2K_{\Ff})$. By \cite[Lemma 3.7]{AD14}, $\Omega^1_{X}|_S\otimes\Omega^1_{X}|_S\rightarrow\mathcal{O}_S(2K_{\Ff})$ extends uniquely to a map
$\Omega^1_{S^\nu}\otimes\Omega^1_{S^\nu}\rightarrow\mathcal{O}_{S^\nu}(2K_{\Ff})$. By construction of $\Ff_S$, $\Omega^1_{S^\nu}\otimes\Omega^1_{S^\nu}\rightarrow\mathcal{O}_{S^\nu}(2K_{\Ff})$ factors through $\mathcal{O}_{S^\nu}(2K_{\Ff_S})$. Thus there exists a Weil divisor $M\geq 0$ on $S^\nu$ such that
$$2K_{\Ff}|_{S^\nu}=2K_{\Ff_S}+M.$$
In particular, near the generic point of $D$, there exists a positive integer $d_0$ such that $M=d_0D$. Then $d_0,\dots,d_m$ satisfy our requirements.
\end{proof}

\subsection{Threefold rank two case}

\begin{thm}[Threefold rank two case]\label{thm: adjunction formula corank 1 dim 3 intro}
Let $(X,\Ff,B=\sum_{j=1}^mb_jB_j)$ be a dlt foliated triple such that $\dim X=3,\rk\Ff=2$, and $B_j$ are the irreducible components of $B$. Assume that
\begin{itemize}
    \item $S$ is an $\Ff$-invariant surface in $X$,
    \item $S^\nu$ is the normalization of $S$,
    \item  $\widehat{X}$ is the formal completion of $X$ along $S$,
    \item $D_1,\dots,D_u$ are all $\Ff$-invariant divisors on $\widehat{X}$ that are not equal to $S$, and
    \item $K_X,S,\sum_{i=1}^uD_i$ are $\Qq$-Cartier, and  $B_j$ is $\Qq$-Cartier for any $j$.
\end{itemize} 
Then there exist prime divisors $C_1,\dots,C_n,T_1,\dots,T_l$ on $S^\nu$, positive integers $w_1,\dots,w_n,\lambda_1,\dots,\lambda_l$, $\{e_{i,k}\}_{1\leq i\leq l,1\leq k\leq n_i}$, and non-negative integers $\{w_{i,j}\}_{1\leq i\leq n,1\leq j\leq m},\{e_{i,k,j}\}_{1\leq i\leq l, 1\leq k\leq n_i,1\leq j\leq m}$, such that for any real numbers $b_1',\dots,b_m'$, we have the following.
\begin{enumerate}
    \item 
    $$\left(K_{\Ff}+\sum_{j=1}^mb_j'B_j\right)\Bigg|_{S^\nu}=K_{S^\nu}+\sum_{i=1}^n\frac{w_i-1+\sum_{j=1}^mw_{i,j}b_j'}{w_i}C_i+\sum_{i=1}^l\lambda_iT_i.$$
    \item $$\left(K_{\widehat{X}}+S+\sum_{i=1}^uD_i+\sum_{j=1}^mb_j'B_j\right)\Bigg|_{S^\nu}=K_{S^\nu}+\sum_{i=1}^n\frac{w_i-1+\sum_{j=1}^mw_{i,j}b_j'}{w_i}C_i+\sum_{i=1}^lT_i.$$
    \item If $(X,\Ff,\sum_{j=1}^mb_j'B_j)$ is lc, then
    $$\left(S^\nu,\sum_{i=1}^n\frac{w_i-1+\sum_{j=1}^mw_{i,j}b_j'}{w_i}C_i+\sum_{i=1}^lT_i\right)$$
    is lc.
    \item For any $i$ such that $\lambda_i\geq 2$,
    \begin{enumerate}
        \item $T_i$ is an lc center of $(X,\Ff,0)$,
        \item $\Ff$ has simple singularities at the generic point of $T_i$, and
        \item for any general closed point $x$ in $T_i$, there are exactly two different separatrices of $\Ff$ containing $x$. One of these two separatrices is a strong separatrix of $\Ff$ at $x$, and the other is $S$, which is not a strong separtrix of $\Ff$ at $x$.
    \end{enumerate}
    \item For any $i$ such that $S$ is a strong separatrix along any general closed point in $T_i$, $\lambda_i=1$. 
    \item Let $T_i^\nu$ be the normalization of $T_i$ for each $i$. For any $i$ such that $\lambda_i\geq 2$, there exist a Weil divisor $P_i\geq 0$ on $T_i^\nu$ and prime divisors (closed points) $P_{i,1},\dots,P_{i,n_i}$ on $T_i^\nu$ satisfying the following.
    \begin{enumerate}
        \item   $$\left(K_{\Ff}+\sum_{j=1}^mb_j'B_j\right)\Bigg|_{T_i^\nu}=K_{T_i^\nu}+P_i+\sum_{k=1}^{n_i}\frac{e_{i,k}-1+\sum_{j=1}^me_{i,k,j}b_j'}{e_{i,k}}P_{i,k}.$$
        \item Any component of $\Supp P_i$ is an lc center of $(X,\Ff,0)$.
        \item If $(X,\Ff,\sum_{j=1}^mb_j'B_j)$ is lc, then 
        $$\left(T_i^\nu,(P_{i})_{\red}+\sum_{k=1}^{n_i}\frac{e_{i,k}-1+\sum_{j=1}^me_{i,k,j}b_j'}{e_{i,k}}P_{i,k}\right)$$ 
        is lc.
    \end{enumerate}
\end{enumerate}
\end{thm}

\begin{proof}
The proof can be done by following the same lines of the proofs of \cite[Lemma 3.18, Corollary 3.20, Lemma 3.22]{CS21} and \cite[Proposition 3.10]{Che22}. Here we provide a short proof by only using the results in \cite{Spi20,CS21} and not their proofs.

\medskip

\noindent\textbf{Step 1}. We prove (2) and (3). By \cite[Lemma 8.14]{Spi20}, $(\widehat{X},S+\sum_{i=1}^uD_i+B)$ is lc. We let
$$K_{S^\nu}+D(\Delta):=\left(K_{\widehat{X}}+S+\sum_{i=1}^uD_i+\Delta\right)\Biggm|_{S^\nu}$$
for any $\Rr$-Cartier $\Rr$-divisor $\Delta$. Let $C$ be a prime divisor on $S^\nu$ and $\eta_C$ the generic point of $C$. If $(\widehat{X},S+\sum_{i=1}^uD_i)$ is lc but not klt at $\eta_C$, then $\eta_C$ is not contained in $B_j$ for any $j$. By \cite[16.7 Corollary]{Kol+92}, the coefficient of $C$ in $D(\Delta)$ is either $0$ or $1$ for any $\Rr$-divisor $\Delta$ such that $\Supp\Delta=\Supp B$. If $(\widehat{X},S+\sum_{i=1}^uD_i)$ is klt at $\eta_C$, then by \cite[Theorem 3.10]{HLS19}, there exist  positive integer $w_C$ and non-negative integers $w_{C,j}$, such that for any real numbers $b_1',\dots,b_m'$, the coefficient of $C$ in $D(\sum_{j=1}^mb_j'B_j)$ is equal to $$\frac{w_C-1+\sum_{j=1}^mw_{C,j}b_j'}{w_C}.$$ 
This implies (2). (3) follows from (2) and \cite[Theorem 3.10]{HLS19} (\cite[16.9 Proposition]{Kol+92} for the $\Qq$-coefficient case).

In the following, we let $C_1,\dots,C_n.T_1,\dots,T_l,w_1,\dots,w_n$, and $\{w_{i,j}\}_{1\leq i\leq n, 1\leq j\leq m}$ be as in (2).

\medskip

\noindent\textbf{Step 2}. We prove (1).

\begin{claim}\label{claim: (1) of adjunction for lc case}
Suppose that $(X,\Ff,\sum_{j=1}^mb_j'B_j)$ is lc. Then (1) holds.
\end{claim}
\begin{proof}
By Theorem \ref{thm: rational polytope foliation intro}, there exist vectors $\bm{b}_1,\dots,\bm{b}_{m+1}\in\mathbb Q^{m}$ and real numbers $a_1,\dots,a_{m+1}\in (0,1]$, such that $\bm{b}_i=(b_{i,1},\dots,b_{i,m})$, $(X,\Ff,B^i:=\sum_{j=1}^mb_{i,j}B_j)$ is lc for any $i$, $\sum_{i=1}^{m+1}a_i=1$, and $\sum_{i=1}^{m+1}a_i\bm{b}_i=(b_1',\dots,b_m')$. By \cite[Corollary 3.20]{CS21}, 
 $$\left(K_{\Ff}+B^k\right)\Bigg|_{S^\nu}=K_{S^\nu}+\sum_{i=1}^n\frac{w_i-1+\sum_{j=1}^mw_{i,j}b_{k,j}}{w_i}C_i+\sum_{i=1}^l\lambda_{i,k}T_i.$$
where $\lambda_{i,k}$ are positive integers.

 Suppose that $\lambda_{r,k}\not=\lambda_{r,k'}$ for some $r$ and some $k\not=k'$. Let $N$ be a sufficiently large positive integer, then by \cite[Corollary 3.20]{CS21} and (2), there exist positive integers $\mu_1,\dots,\mu_{l}$, such that
 \begin{align*}
    &K_{S^\nu}+\sum_{i=1}^n\frac{w_i-1+\sum_{j=1}^mw_{i,j}(\frac{1}{N}b_{k,j}+\frac{N-1}{N}b_{k',j})}{w_i}C_i+\sum_{i=1}^l\mu_iT_i=\left(K_{\Ff}+\frac{1}{N}B^k+\frac{N-1}{N}B^{k'}\right)\Bigg|_{S^\nu}\\
    =&K_{S^\nu}+\sum_{i=1}^n\frac{w_i-1+\sum_{j=1}^mw_{i,j}(\frac{1}{N}b_{k,j}+\frac{N-1}{N}b_{k',j})}{w_i}C_i+\sum_{i=1}^l\frac{1}{N}(\lambda_{i,k}+(N-1)\lambda_{i,k'})T_i,
 \end{align*}
 which is not possible as $\mu_r$ is an integer but $\frac{1}{N}(\lambda_{r,k}+(N-1)\lambda_{r,k'})$ is not. Thus for any $i$, $\lambda_{i,k}$ is a constant for any $k$, so we may let $\lambda_i:=\lambda_{i,k}$ for any $i$. (1) immediately follows in this case.
\end{proof} 

\noindent\textit{Proof of Theorem \ref{thm: adjunction formula corank 1 dim 3 intro} continued}. (1) follows immediately from Claim \ref{claim: (1) of adjunction for lc case} and linearity of the coefficients of $B_1,\dots,B_m$. 

\medskip

\noindent\textbf{Step 3}. We conclude the proof in this step.

(4) By  \cite[Corollary 3.20]{CS21}, the generic point of $T_i$ is not contained in $\Sing(X)$.

Suppose that $T_i$ is not an lc center of $(X,\Ff,0)$. Then $(X,\Ff,0)$ is terminal at the generic point of $T_i$. By \cite[Lemma 3.14]{CS21}, $S$ is the unique $\Ff$-invariant divisor passing through the generic point of $T_i$. Thus $\sum_{i=1}^uD_i=0$ near the generic point of $T_i$. By (2), the generic point of $T_i$ is contained in $\Sing(X)$, a contradiction. Thus $T_i$ is an lc center of $(X,\Ff,0)$, which implies (4.a). Since $(X,\Ff,0)$ is dlt, $\Ff$ has simple singularities near the generic point of $T_i$, which is (4.b). (4.c) follows from (4.b), \cite[Lemma 3.14]{CS20}, and  \cite[Corollary 3.20]{CS21}.

(5) It immediately follows from (4).

(6) By (4), we may let $S'$ be the strong separatrix along $T_i$ and $S'^\nu$ the normalization of $S'$. Then by (1)(3)(5),
$$\left(K_{\Ff}+\sum_{j=1}^mb_j'B_j\right)\Bigg|_{S'^\nu}=K_{S'^\nu}+\sum_k\frac{w'_k-1+\sum_{j=1}^mw'_{k,j}b_j'}{w_k'}C_k'+\sum_k\lambda'_kT_k'.$$
for some prime divisors $C_k',T_k'$, positive integers $w_k',\lambda_k'$, and non-negative integers $w_{k,j}'$, such that the images of $T_i'$ and $T_i$ in $X$ coincide, $\lambda_i'=1$, and
$$\left(S'^\nu,\sum_k\frac{w'_k-1+\sum_{j=1}^mw'_{k,j}b_j'}{w_k'}C_k'+\sum_kT_k'\right)$$
is lc. In particular, the normalization of $T_i'$ is $T_i^\nu$. By \cite[16.6.3, 16.7 Corollary]{Kol+92} and \cite[Theorem 3.10]{HLS19}, we have
\begin{align*}
    (K_{\Ff}+\sum_{j=1}^mb_j'B_j)|_{T_i^\nu}&=(K_{S'^\nu}+\sum_k\frac{w'_k-1+\sum_{j=1}^mw'_{k,j}b_j'}{w'_k}C_k'+\sum_k\lambda'_kT_k')|_{T_i^\nu}\\
    &=K_{T_i^\nu}+P_i'+\sum_{k=1}^{n_i}\frac{e_{i,k}-1+\sum_{j=1}^me_{i,k,j}b_j'}{e_{i,k}}P_{i,k}+\sum_{\lambda_k'\geq 2}(\lambda_k'-1)T_k'|_{T_i^\nu}.
\end{align*}
for some non-negative integer $n_i$, positive integers $\{e_{i,k}\}_{k=1}^{n_i},\{e_{i,k,j}\}_{1\leq k\leq n_i,1\leq j\leq m}$, and an effective Weil divisor $P_i\geq 0$, such that
\begin{itemize}
\item $$\left(T_i^\nu,P_i'+\sum_{k=1}^{n_i}\frac{e_{i,k}-1+\sum_{j=1}^me_{i,k,j}b_j'}{e_{i,k}}P_{i,k}\right)$$ is lc, and
\item $\Supp T_k'|_{T_i^\nu}\subset\Supp P_i'$ for any $k\not=i$.
\end{itemize}
For any point $D\subset\Supp P_i'$, by \cite[Lemma 3.22]{CS21} (applied to the case when $b_j'=0$ for every $j$), $D$ is an lc center of $(X,\Ff,0)$. Since $(X,\Ff,0)$ is dlt, $T_i$ is Cartier near $D$. Thus $\sum_{\lambda_k'\geq 2}(\lambda_k'-1)T_k'|_{T_i^\nu}$
is an effective Weil divisor and we may let $P_i:=P_i'+\sum_{\lambda_k'\geq 2}(\lambda_k'-1)T_k'|_{T_i^\nu}$.
(6.a) and (6.b) immediately follow, and (6.c) follows from the adjunction formula for usual pairs.
\end{proof}

\section{Uniform rational polytopes}\label{sec: Uniform rational polytopes}

In this section, we establish the theory of uniform rational polytopes and functional divisors for foliations in dimension $\leq 3$ and prove Theorem \ref{thm: uniform rational polytope foliation intro}. 

\subsection{Log canonical thresholds}\label{subsec: lct}

In this subsection, we define the lc thresholds of foliations. Our definition is slightly different from the traditional definition since we need to consider lc thresholds of $\Rr$-divisors which are not necessarily effective.

\begin{defn}[Lc thresholds]
Let $(X,\Ff,B)$ be an lc foliated sub-triple and $D$ an $\Rr$-Cartier $\Rr$-divisor on $X$. An \emph{lc threshold} (\emph{lct} for short) of $D$ with respect to $(X,\Ff,B)$ is a real number $t_0$, such that
\begin{enumerate}
    \item $(X,\Ff,B+t_0D)$ is sub-lc, and
    \item for any positive real number $\delta$, either $(X,\Ff,B+(t_0+\delta)D)$ or  $(X,\Ff,B+(t_0-\delta)D)$ is not sub-lc.
\end{enumerate}
When $D\geq 0$, the lc threshold of $D$ with respect to $(X,\Ff,B)$ is unique, and we denote the lc threshold of $D$ with respect to $(X,\Ff,B)$ by $\lct(X,\Ff,B;D)$. 
\end{defn}

The following lemma indicates that lc thresholds can always be achieved in dimension $\leq 3$.

\begin{lem}\label{lem: lc threshold has lc center}
Let $(X,\Ff,B)$ be a sub-lc foliated sub-triple of dimension $\leq 3$ and $D$ an $\Rr$-Cartier $\Rr$-divisor on $X$. Let $B(t):=B+tD$ for any real number $t$. Then there exist $t_1,t_2\in\mathbb R\cup\{-\infty,+\infty\}$, such that
\begin{enumerate}
\item $t_1\leq 0\leq t_2$,
\item for any real number $t$, $(X,\Ff,B(t))$ is sub-lc if and only if $t_1\leq t\leq t_2$, and
\item if $t_1\not=t_2$, then for any $i\in\{1,2\}$,
\begin{itemize}
    \item either $t_i\in\{-\infty,+\infty\}$, or
    \item there exists a prime divisor $E_i$ over $X$, such that $a(E_i,X,\Ff,B(t_i))=-\epsilon_{\Ff}(E_i)$ and $a(E_i,X,\Ff,B(t))>-\epsilon_{\Ff}(E_i)$ for any $t\in (t_1,t_2)$.
\end{itemize} 
\end{enumerate}
\end{lem}
\begin{proof}
We may let $t_1:=\inf\{t\mid (X,\Ff,B(t))\text{ is sub-lc}\}$ and $t_2:=\sup\{t\mid (X,\Ff,B(t))\text{ is sub-lc}\}$. (1) immediately follows.

By \cite[Theorem 4.5]{LLM23}, there exists a foliated log resolution $f: Y\rightarrow X$ of $(X,\Ff,B)$. Let $K_{\Ff_Y}+B_Y(t):=f^*(K_{\Ff}+B(t))$ for any real number $t$ where $\Ff_Y:=f^{-1}\Ff$. Possibly replacing $(X,\Ff,B(t))$ with $(Y,\Ff_Y,B_Y(t))$, we may assume that $(X,\Ff,\Supp B\cup\Supp D)$ is foliated log smooth. By \cite[Lemma 4.3]{LLM23}, $(X,\Ff,B(t))$ is lc if and only if for any component $T$ of $\Supp B\cup\Supp D$, $\mult_TB(t)\leq\epsilon_{\Ff}(T)$.
Since the coefficients of $B(t)$ are affine functions, $t_1=\min\{t\mid (X,\Ff,B(t))\text{ is sub-lc}\}$ or $-\infty$, and $t_2=\max\{t\mid (X,\Ff,B(t))\text{ is sub-lc}\}$ or $+\infty$. This implies (2). (3) immediately follows from (2).
\end{proof}

\subsection{Special dlt models}\label{subsec: special dlt}
In this subsection we prove some key lemmas related to dlt models of different foliated triples.

\begin{lem}\label{lem: dlt not influence perturbation}
Let $c,m$ be be positive integers, $r_1,\dots,r_c$ real numbers such that $1,r_1,\dots,r_c$ are linearly independent over $\mathbb Q$, $\bm{r}:=(r_1,\dots,r_c)$, and $s_1,\dots,s_m: \mathbb R^{c+1}\rightarrow\mathbb R$ $\mathbb Q$-linear functions. Assume that
\begin{itemize}
    \item $(X,\Ff,B(\bm{r}):=\sum_{i=1}^ms_i(1,\bm{r})B_i)$ is an lc foliated triple of dimension $\leq 3$, and
    \item $B_i\geq 0$ are distinct Weil divisors and $s_i(1,\bm{r})\geq 0$ for each $i$.
\end{itemize} 
Let $f: Y\rightarrow X$ be a dlt modification of $(X,\Ff,B(\bm{r}))$ with prime exceptional divisors $E_1,\dots,E_n$, $\Ff_Y:=f^{-1}\Ff$, $B(\bm{v}):=\sum_{i=1}^ms_i(1,\bm{v})B_i$ for any $\bm{v}\in\mathbb R^c$, and $B_Y(\bm{v}):=f^{-1}_*B(\bm{v})+\sum_{i=1}^n\epsilon_{\Ff}(E_i)E_i$. Then for any $\bm{v}\in\mathbb R^c$,
\begin{enumerate}
  \item $K_{\Ff_Y}+B_Y(\bm{v})=f^*(K_{\Ff}+B_Y(\bm{v}))$, and
    \item $(X,\Ff,B(\bm{v}))$ is lc if and only if  $(Y,\Ff_Y,B_Y(\bm{v}))$ is lc.
\end{enumerate}
\end{lem}
\begin{proof}
There exist $\Qq$-divisors $M_0,\dots,M_c$, such that $B(1,v_1,\dots,v_c)=M_0+\sum_{i=1}^cv_iM_i$ for any $v_1,\dots,v_c$. By \cite[Lemma 5.3]{HLS19}, $M_i$ is $\Qq$-Cartier for any $1\leq i\leq m$. Then for any $j$,
$$-\epsilon_{\Ff}(E_j)=a(E_i,\Ff,B(\bm{r}))=a(E_j,\Ff,M_0)-\sum_{i=1}^mr_i\mult_{E_j}M_i,$$
so $\mult_{E_j}M_i=0$ for each $i$. This implies (1) and (2) follows from (1).
\end{proof}

\begin{lem}\label{lem: find dlt model with nontrivial multiplicity}
Let $(X,\Ff,B)$ be an lc foliated triple and $M$ an $\Rr$-Cartier $\Rr$-divisor on $X$, such that for any positive real number $\delta$, either $(X,\Ff,B+\delta M)$ is not lc or $(X,\Ff,B-\delta M)$ is not lc. Assume that
\begin{itemize}
    \item either $\dim X\leq 3$ and $\rk\Ff=\dim X-1$, or
    \item $\dim X=3$, $\rk\Ff=1$, and $X$ is projective.
\end{itemize}
Then one of the following holds.
\begin{enumerate}
    \item There exists a prime divisor $D$ on $X$, such that $a(D,\Ff,B)=-\epsilon_{\Ff}(D)$ and $\mult_DM\not=0$.
    \item There exists a dlt modification $f: Y\rightarrow X$ of $(X,\Ff,B)$ with prime $f$-exceptional divisors $E_1,\dots,E_n$ satisfying the following.
    \begin{enumerate}
        \item Let $K_{\Ff_Y}+B_Y:=f^*(K_{\Ff}+B)$, where $\Ff_Y:=f^{-1}\Ff$. Then there exists a positive real number $\delta$ such that $(Y,\Ff_Y,B_Y+\delta\{B_Y\})$ is dlt.
         \item There exists an integer $1\leq i\leq n$, such that $\mult_{E_i}M\not=0$.
    \end{enumerate}
\end{enumerate}
\end{lem}
\begin{proof}
By Lemma \ref{lem: lc threshold has lc center}(2), there exists a prime divisor $E$ that is exceptional over $X$, such that $a(E,X,\Ff,B)=-\epsilon_{\Ff}(E)$ and $\mult_EM\not=0$. If $E$ is on $X$, then either (1) holds, so we may assume that $E$ is exceptional over $X$. By \cite[Theorem 4.5]{LLM23}, we may let $g: Z\rightarrow X$ be a log resolution of $(X,\Ff,\Supp B\cup\Supp M)$ such that $E$ is on $Z$. Let $F_1,\dots,F_m$ be the prime $g$-exceptional divisors, $\Ff_Z:=g^{-1}\Ff$,  $B_Z:=g^{-1}_*B$, and $M_Z:=g^{-1}_*M$. By \cite[Lemma 4.3]{LLM23} and the definition of dlt singularities, $$\left(Z,\Ff_Z,\Supp(B_Z+M_Z)+\sum_{i=1}^m\epsilon_{\Ff}(F_i)\right)$$
is dlt, and
$$K_{\Ff_Z}+B_Z-\sum_{i=1}^ma(F_i,\Ff,B)F_i=g^*(K_{\Ff}+B).$$
Let $F:=\sum_{i=1}^m(\epsilon_{\Ff}(F_i)+a(F_i,\Ff,B))F_i$. Then $F\geq 0$ is exceptional over $X$ and $E\not\subset\Supp F$. Since $K_{\Ff_Z}+B_Z+\sum_{i=1}^m\epsilon_{\Ff}(F_i)F_i\sim_{\mathbb R,X}F$, by \cite[Theorem 5.9]{LLM23}, we may run a $(K_{\Ff_Z}+B_Z+\sum_{i=1}^m\epsilon_{\Ff}(F_i)F_i)$-MMP$/X$ which terminates with a good minimal model (cf. \cite[Definition 5.5]{LLM23}) $(Y,\Ff_Y,B_Y)/X$ of $(Z,\Ff_Z,B_Z+\sum_{i=1}^m\epsilon_{\Ff}(F_i)F_i)/X$ such that $K_{\Ff_Y}+B_Y\sim_{\mathbb R,X}0$ and the divisors contracted by $Z\dashrightarrow Y$ are exactly those contained in $\Supp F$. Then $(Y,\Ff_Y,B_Y)$ is $\Qq$-factorial dlt and the induced birational morphism $f: Y\rightarrow X$ is a dlt modification of $(X,\Ff,B)$.  

Since $(Z,\Ff_Z,\Supp(B_Z+M_Z)+\sum_{i=1}^m\epsilon_{\Ff}(F_i))$ is dlt, there exists a positive real number $\delta$ such that $(Z,\Ff_Z,B_Z+\sum_{i=1}^m\epsilon_{\Ff}(F_i)F_i+\delta\{B_Z\})$ is dlt and the induced birational map $Z\dashrightarrow Y$ is also a partial $(K_{\Ff_Z}+B_Z+\sum_{i=1}^m\epsilon_{\Ff}(F_i)F_i+\delta\{B_Z\})$-MMP for some positive real number $\delta$. By \cite[Theorem 5.8]{LLM23}, $(Y,\Ff_Y,B_Y+\delta\{B_Y\})$ is dlt. Thus (2.a) holds.

Let $E_1,\dots,E_n$ be the prime $f$-exceptional divisors. Since $E\not\subset\Supp F$, we may let $i$ be the index such that $E_i$ is the image of $E$ on $Y$. Thus (2.b) holds.
\end{proof}

\begin{lem}[Key Lemma]\label{lem: find nontrivial divisor on dlt model}
Let $(X,\Ff,B)$ be an lc foliated triple and $M$ an $\Rr$-Cartier $\Rr$-divisor on $X$, such that 
\begin{itemize}
    \item 
    \begin{itemize}
        \item either $\dim X\leq 3$ and $\rk\Ff=\dim X-1$, or
        \item  $\dim X=3$, $\rk\Ff=1$, and $X$ is projective,
    \end{itemize} 
    \item $(X,\Ff,B+M)$ is lc,
    \item $(X,\Ff,B+(1+\epsilon)M)$ is not lc for any positive real number $\epsilon$,
    \item $\Supp B=\Supp(B+M)$, and
    \item for any prime divisor $D$ on $X$ such that $a(D,\Ff,B+M)=-\epsilon_{\Ff}(D)$, $\mult_DM=0$.
\end{itemize}
Then there are two projective birational morphisms $h: X'\rightarrow X$ and $g: Y'\rightarrow X'$ and a real number $t\in (0,1)$ satisfying the following. 
\begin{enumerate}
    \item $h$ is a dlt modification of $(X,\Ff,B+tM)$.
    \item For any prime $h$-exceptional divisor $D$, $a(D,\Ff,B)=-\epsilon_{\Ff}(D)$. In particular, $\mult_DM=0$ and $a(D,\Ff,B+sM)=-\epsilon_{\Ff}(D)$ for any real number $s$.
    \item $g$ extracts a unique prime divisor $E$. In particular, $-E$ is ample over $X'$.
    \item $a(E,\Ff,B+M)=-\epsilon_{\Ff}(E)$ and $a(E,\Ff,B)>-\epsilon_{\Ff}(E)$. In particular, $\mult_EM>0$ and $a(E,\Ff,B+sM)>-\epsilon_{\Ff}(E)$ for any real number $s<1$.
    \item Let $B_{Y'},M_{Y'}$ be the strict transforms of $B,M$ on $Y'$ respectively, $\Ff_{Y'}:=(h\circ g)^{-1}\Ff$, and
    $$F_{Y'}:=\sum_{D\text{ is a prime }h\circ g\text{-exceptional divisor}}\epsilon_{\Ff}(D)D.$$
    Then $(Y',\Ff_{Y'},B_{Y'}+tM_{Y'}+F_{Y'})$ is $\Qq$-factorial dlt.
\end{enumerate}
\end{lem}
\begin{proof}
By Lemma \ref{lem: find dlt model with nontrivial multiplicity}, there exists a dlt modification $f: Y\rightarrow X$ of $(X,\Ff,B+M)$ satisfying the following: let $\Ff_Y:=f^{-1}\Ff$, $B_Y:=f^{-1}_*B,M_Y:=f^{-1}_*M$, $K_{\Ff_Y}+B_Y:=f^*(K_{\Ff}+B)$, and $E_1,\dots,E_n$ the prime $f$-exceptional divisors, then
\begin{itemize}
    \item $(Y,\Ff_Y,B_Y+M_Y+\delta\{B_Y+M_Y\}+\sum_{i=1}^n\epsilon_{\Ff}(E_i)E_i)$ is dlt for some positive real number $\delta$, and
    \item $\mult_{E_i}M\not=0$ for some $1\leq i\leq n$.
\end{itemize}
Since $\Supp B=\Supp(B+M)$, $(Y,\Ff_Y,B_Y+tM_Y+\sum_{i=1}^n\epsilon_{\Ff}(E_i)E_i)$ is dlt for some $t\in (0,1)$. Since $(X,\Ff,B)$ and $(X,\Ff,B+M)$ are lc, $(X,\Ff,B+tM)$ is lc. Thus
$$K_{\Ff_Y}+B_Y+tM_Y+\sum_{i=1}^n\epsilon_{\Ff}(E_i)E_i\sim_{\mathbb R,X}\sum_{i=1}^n(\epsilon_{\Ff}(E_i)+a(E_i,\Ff,B_Y+tM_Y))E_i\geq 0.$$
By \cite[Theorem 5.9]{LLM23}, we may run a $(K_{\Ff_Y}+B_Y+tM_Y+\sum_{i=1}^n\epsilon_{\Ff}(E_i)E_i)$-MMP$/X$ which terminates with a good minimal model $(X',\Ff',B'+tM'+F')/X$ of $(Y,\Ff_Y,B_Y+tM_Y+\sum_{i=1}^n\epsilon_{\Ff}(E_i)E_i)/X$ such that $K_{\Ff'}+B'+tM'+F'\sim_{\mathbb R,X}0$, where $B',M',F'$ are the images of $B_Y,M_Y,\sum_{i=1}^n\epsilon_{\Ff}(E_i)E_i$ on $X'$ respectively. Then $(X',\Ff',B'+tM'+F')$ is $\Qq$-factorial dlt and the induced morphism $h: X'\rightarrow X$ is a dlt modification of $(X,\Ff,B+tM)$. 

By construction, the divisors contracted by the induced birational map $Y\dashrightarrow X'$ are all divisors $E_i$ such that $\epsilon_{\Ff}(E_i)>a(E_i,\Ff,B_Y+tM_Y)$. Since $\mult_{E_i}M\not=0$, $Y\dashrightarrow X'$ contracts $E_i$, hence $Y\dashrightarrow X'$ contains a divisorial contraction. We let $g: Y'\dashrightarrow X'$ be the last step of the $(K_{\Ff_Y}+B_Y+tM_Y+\sum_{i=1}^n\epsilon_{\Ff}(E_i)E_i)$-MMP$/X$. Since $X'$ is $\Qq$-factorial and $K_{\Ff'}+B'+tM'+F'\sim_{\mathbb R,X}0$, $g$ is a divisorial contraction of a prime divisor $E$. 

We show that $h,g$ and $t$ satisfy our requirements.

(1)(5) immediately follow from our construction. 

For any prime divisor $D$ on $X'$ that is exceptional over $X$, the center of $D$ on $Y$ is a divisor that is exceptional over $X$. Since $f: Y\rightarrow X$ is a dlt modification of $(X,\Ff,B+M)$, $a(D,\Ff,B+M)=-\epsilon_{\Ff}(D)$. By (1), $h$ is a dlt modification of $(X,\Ff,B+tM)$, so  $a(D,\Ff,B+tM)=-\epsilon_{\Ff}(D)$. Thus $\mult_DM=0$, so $a(D,\Ff,B)=-\epsilon_{\Ff}(D)$. This implies (2). 

Since $g$ is a divisorial contraction of a prime divisor $E$, by the negativity lemma, $-E$ is ample$/X'$. This implies (3). 

Since the center of $E$ on $Y$ is a divisor that is exceptional over $X$, $a(E,\Ff,B+M)=-\epsilon_{\Ff}(E)$. Since $E$ is contracted by the $(K_{\Ff_Y}+B_Y+tM_Y+\sum_{i=1}^n\epsilon_{\Ff}(E_i)E_i)$-MMP$/X$: $Y\dashrightarrow X'$, $$a(E,\Ff,B+tM)=a(E,\Ff',B'+tM'+F')>-\epsilon_{\Ff}(E).$$
Thus $a(E,\Ff,B)>-\epsilon_{\Ff}(E)$. This implies (4).
\end{proof}

\subsection{Threefold rank one case}\label{subsec: 3fold rk1}

In this subsection, we prove a slightly weaker version of Theorem \ref{thm: uniform rational polytope foliation intro} when $\dim X=3$ and $\rk\Ff=1$.

\begin{prop}\label{prop: uniform rational polytope foliation one variable threefold rk1}
Let $c,m$ be positive integers, $r_1,\dots,r_c$ real numbers such that $1,r_1,\dots,r_c$ are linearly independent over $\mathbb Q$, $\bm{r}:=(r_1,\dots,r_c)$, and $s_1,\dots,s_m: \mathbb R^{c+1}\rightarrow\mathbb R$ $\mathbb Q$-linear functions. Then there exists a positive real number $\delta$ depending only on $\bm{r}$ and $s_1,\dots,s_m$ satisfying the following. Assume that
\begin{enumerate}
    \item $(X,\Ff,B=\sum_{i=1}^ms_i(1,\bm{r})B_i)$ is a projective lc foliated triple such that $\dim X=3$ and $\rk\Ff=1$,
    \item $B_i\geq 0$ are distinct Weil divisors (possibly $0$) and $s_i(1,\bm{r})\geq 0$, and
    \item $B(t):=\sum_{i=1}^ms_i(1,r_1,\dots,r_{c-1},t)B_i$ for any $t\in\mathbb R$.
\end{enumerate}
Then $(X,\Ff,B(t))$ is lc for any $t\in (r_c-\delta,r_c+\delta)$.
\end{prop}
\begin{proof}
We let $s_i(t):=s_i(1,r_1,\dots,r_{c-1},t)$ for any $t\in\mathbb R$. If $s_i(r_c)=0$, then $s_i(t)=0$ for any $t$, so we may assume that $s_i(r_c)\not=0$ for any $i$. By Lemma \ref{lem: dlt not influence perturbation} and \cite[Theorem 5.7]{LLM23}, possibly replacing $(X,\Ff,B(r_c))$ with a dlt model, we may assume that $(X,\Ff,B(r_c))$ is $\Qq$-factorial dlt. In particular, $\Ff$ is non-dicritical.

We only need to prove that there exists a positive real number $\epsilon$ depending only on $\bm{r}$ and $s_1,\dots,s_m$, such that for any lc threshold $t_0$ of $(X,\Ff,B(t))$, $|t_0-r_c|>\epsilon$. Thus we may assume that $(X,\Ff,B(t))$ has an lc threshold $t_0$. Since $1,r_1,\dots,r_c$ are linearly independent over $\Qq$, $r_c\not=t_0$. Moreover, there exists a positive real number $\delta_1$ depending only on $\bm{r}$ and $s_1,\dots,s_m$, such that for any $t\in (r_c-\delta_1,r_c+\delta_1)$ and any $i$, $s_i(t)\geq\frac{1}{2}s_i(r_c)>0$. In particular, for any $t\in (r_c-\delta_1,r_c+\delta_1)$, $\Supp B(t)=\Supp B(r_c)$, $\Supp\lfloor B(t)\rfloor=\Supp\lfloor B(r_c)\rfloor$, and $B(t)\geq\frac{1}{2}B(r_c)$.

We may assume that $t_0\in (r_c-\delta_1,r_c+\delta_1)$. In particular, we may assume that $(X,\Ff,B(t_0))$ does not have an lc center in codimension $1$ that is also an lc center of $(X,\Ff,B(r_c))$. Thus $(X,\Ff,B(t_0))$ has an lc center $x$ such that $\dim\bar x\leq 1$, and $x$ is not an lc center of $(X,\Ff,B(r_c))$. In particular, $x\in\Supp B(r_c)$. 

By Lemma \ref{lem: find nontrivial divisor on dlt model} (more precisely, $X,\Ff,B,M$ in Lemma \ref{lem: find nontrivial divisor on dlt model} correspond to our $X,\Ff,B(r_c),B(t_0)-B(r_c)$ respectively), possibly replacing $(X,\Ff,B(t))$, we may assume that there exists a divisorial contraction $g: Y\rightarrow X$ of a prime divisor $E$ and a real number $s$ satisfying the following. Let $B_Y(t)$ be the strict transform of $B(t)$ on $Y$ for any $t$ and $\Ff_Y:=g^{-1}\Ff$, then:
\begin{itemize}
\item  $s\in (r_c,t_0)$ if $r_c<t_0$, and $s\in (t_0,r_c)$ if $t_0<r_c$,
\item $(X,\Ff,B(s))$ is $\Qq$-factorial dlt, $(X,\Ff,B(r_c))$ is lc, and  $(X,\Ff,B(t_0))$ is lc. In particular, since $\rk\Ff=1$, by definition, $(X,\Ff,B(r_c))$ and $(X,\Ff,B(t_0))$ are dlt.
\item $-E$ is ample over $X$,
\item $(Y,\Ff_Y,B_Y(s)+\epsilon_{\Ff}(E)E)$ is $\Qq$-factorial dlt, and
\item $a(E,\Ff,B(t_0))=-\epsilon_{\Ff}(E)$ and $a(E,\Ff,B(r_c))>-\epsilon_{\Ff}(E)$. In particular,  $(Y,\Ff_Y,B_Y(t_0)+\epsilon_{\Ff}(E)E)$ is lc.
\end{itemize}

Since $(X,\Ff,B(s))$ is $\Qq$-factorial dlt, $\Ff$ is non-dicritical, so $E$ is $\Ff_Y$-invariant and $\epsilon_{\Ff}(E)=0$. Let $V$ be the normalization of $\Center_XE$, $E^\nu$ the normalization of $E$,  $g|_{E^\nu}: E^\nu\rightarrow V$ the induced projective surjective morphism with Stein factorization
$$E^\nu\xrightarrow{\pi} W\xrightarrow{\tau} V,$$
and $F$ a general fiber of $\pi$. Then since $-E$ is ample over $X$, 
$$K_{\Ff_Y}+B_Y(r_c)=g^*(K_{\Ff}+B(r_c))+a(E,\Ff,B(r_c))E\sim_{\mathbb R,X}a(E,\Ff,B(r_c))E$$
is anti-ample$/X$. Thus
$$(K_{\Ff_Y}+B_Y(r_c))|_F$$
is anti-ample. Moreover, since
$$K_{\Ff_Y}+B_Y(t_0)=g^*(K_{\Ff}+B(t_0))\sim_{\mathbb R,X}0,$$
we have that
$$(K_{\Ff_Y}+B_Y(t_0))|_F\sim_{\mathbb R}0.$$

We let $x$ be the generic point of $\Center_XE$ and let $E^\nu$ be the normalization of $E$. Since $(Y,\Ff_Y,B_Y(s))$ is $\Qq$-factorial dlt, by Theorem \ref{thm: adjunction threefold rank 1 intro}, there exist prime divisors $C_1,\dots,C_n$ on $E^\nu$, non-negative integers $\{w_{i,j}\}_{0\leq i\leq n,1\leq j\leq m}$, and a rank $1$ foliation $\Ff_E$ on $E^\nu$, such that
$$(K_{\Ff}+B_Y(t))|_{E^\nu}=K_{\Ff_E}+\sum_{i=1}^n\frac{w_{i,0}+\sum_{j=1}^mw_{i,j}s_j(t)}{2}C_i$$
for any real number $t$. There are two cases.

\medskip

\noindent\textbf{Case 1}. $\dim\bar x=1$. We let $$h(t):=\sum_{i=1}^n\frac{w_{i,0}+\sum_{j=1}^mw_{i,j}s_j(t)}{2}(C_i\cdot F)$$
for any real number $t$. Then
$K_{\Ff_E}\cdot F+h(t)\not=0$
when $t\not=t_0$ and $K_{\Ff_E}\cdot F+h(t_0)=0$. In particular, $h(t)$ is not a constant function. Since $F$ is an irreducible component of a general fiber, $C_i\cdot F\geq 0$ for any $i$ and $F^2=0$. In particular, $h(t_0)\geq 0$, and $h(t_0)=0$ if and only if $C_i\cdot F=0$ for any $i$. Recall that $s_j(t_0)\geq\frac{1}{2}s_j(r_c)>0$ for any $j$. Therefore, since $h(t)$ is not a constant function, $h(t_0)>0$. There are two cases.

\medskip

\noindent\textbf{Case 1.1}. $F$ is not $\Ff_E$-invariant. In this case, since $E^\nu$ is smooth near $F$,
$$0=K_{\Ff_E}\cdot F+h(t_0)>K_{\Ff_E}\cdot F=(K_{\Ff_E}+F)\cdot F=\tang(\Ff_E,F)\geq 0,$$
a contradiction.

\medskip

\noindent\textbf{Case 1.2}. $F$ is $\Ff_E$-invariant. In this case, since $E^\nu$ is smooth near $F$,
$$0=K_{\Ff_E}\cdot F+h(t_0)>K_{\Ff_E}\cdot F=Z(\Ff_E,F)-\chi(F),$$
so $K_{\Ff_E}\cdot F=Z(\Ff_E,F)-\chi(F)\in\{-1,-2\}$. Thus
$$h(t_0)=\sum_{i=1}^n\frac{c_{i,0}+\sum_{j=1}^mc_{i,j}s_j(t_0)}{2}(C_i\cdot F)\in\{1,2\}.$$
Since $s_j(t_0)\geq\frac{1}{2}s_j(r_c)>0$ for any $j$, there are finitely many possibilities of $c_{i,j}$ (which do not depend on $t_0$) for any $i$ such that $C_i\cdot F\not=0$. Thus there are only finitely many possibilities of $h(t)$. Since $h(t_0)=0$ and $h$ is not a constant function, there are only finitely many possibilities of $t_0$. The proposition follows in this case.

\medskip

\noindent\textbf{Case 2}. $\dim\bar x=0$. In this case $F=E^\nu$ and $x$ is a closed point. We let $$C(t):=\sum_{i=1}^n\frac{w_{i,0}+\sum_{j=1}^mw_{i,j}s_j(t)}{2}C_i.$$ Then
$K_{\Ff_E}+C(t_0)\sim_{\mathbb R}0$, and
$K_{\Ff_E}+C(r_c)$ is anti-ample. Therefore, $C(t_0)\not=0$, so $K_{\Ff_E}$ is not pseudo-effective. By \cite[Theorem 1.1]{CP19}, \cite[Theorem 3.1]{LLM23}, $\Ff_E$ is algebraically integrable. By \cite[3.5]{Dru21}, there exists a projective birational morphism $f: V\rightarrow E^\nu$ and a contraction $\tau: V\rightarrow Z$ to a curve $Z$, such that $\Ff_V:=f^{-1}\Ff_E$ is induced by $\tau$, i.e. $\tau$ is the family of leaves of $\Ff$. Let $L_V$ be a general fiber of $\tau$, $L:=f_*L_V$, and $L_Y$ the image of $L$ on $Y$. 
\begin{claim}\label{claim: 2l cartier}
For any $i$ such that $w_{i,j}\not=0$ for some $j$, $C_i\cdot L$ is well-defined and is an integer.
\end{claim}
\begin{proof}
We only need to show that $L$ is Cartier near any closed point $y\in C_i\cap L$. Since $(Y,\Ff_Y,B_Y(r_c))$ is lc and $s_j(r_c)>0$ for any $j$, $\Supp B_Y(t)$ does not contain any lc center of $(Y,\Ff_Y,0)$ for any real number $t$. Since $w_{i,j}\not=0$ for some $j$, $C_i$ does not contain any lc center of $(Y,\Ff_Y,0)$. In particular, the image of $y$ in $Y$ is not an lc center of $(Y,\Ff_Y,0)$. Since $\Ff_Y$ has simple singularities, $\Ff_Y$ is terminal near $y$. By \cite[Lemma 2.12(1)]{CS20}, $\Ff_Y$ is induced by a fibration near $y$ up to a quasi-\'etale cover. Therefore, there are only finitely many $\Ff_Y$-invariant curves passing through the image of $y$ in $Y$. Thus, there are only finitely many $\Ff_E$-invariant curves passing through $y$, so $y$ is a non-dicritical singularity of $\Ff_E$. Since $L_V$ is a general fiber of $\tau$,  $f$ is an isomorphism near a neighborhood of $y$, and  $L_V$ is smooth near $f^{-1}(y)$. Thus $L$ is smooth near $y$. Therefore,  $L$ is Cartier near $y$ and we are done.
\end{proof}
\noindent\textit{Proof of Proposition \ref{prop: uniform rational polytope foliation one variable threefold rk1} continued}.
We let
$$h(t):=\sum_{i=1}^n\frac{\sum_{j=1}^mw_{i,j}s_j(t)}{2}(C_i\cdot L)$$ for any real number $t$, then $h(r_c)\geq 0$, $h(t_0)\geq 0$, and $h(t)=B_Y(t)\cdot L_Y$ for any $t$. Since
$$0>(K_{\Ff_Y}+B_Y(r_c))\cdot L_Y=K_{\Ff_Y}\cdot L_Y+h(r_c),$$
we have that 
$$\left(K_{\Ff_E}+\sum_{i=1}^n\frac{w_{i,0}}{2}C_i\right)\cdot L=K_{\Ff_Y}\cdot L_Y<0.$$
Since $\Ff_Y$ has simple singularities, for any closed point $y\in L_Y$ such that $\Ff_Y$ is not terminal at $y$, $2K_{\Ff_Y}$ is Cartier near $y$. By \cite[Proposition 2.16(3)(4), Proposition 3.3]{CS20},
$$K_{\Ff_Y}\cdot L_Y=-2+\frac{1}{2}\lambda+\sum_{i=1}^u\frac{\mu_i-1}{\mu_i}$$
where $\lambda,u$ are non-negative integers and $\mu_i$ are positive integers. Therefore,
$$0>(K_{\Ff_Y}+B_Y(r_c))\cdot L_Y=-2+\frac{1}{2}\lambda+\sum_{i=1}^u\frac{\mu_i-1}{\mu_i}+\sum_{i=1}^n\sum_{j=1}^m\frac{w_{i,j}(C_i\cdot L)}{2}s_j(r_c)$$
and
$$0=(K_{\Ff_Y}+B_Y(t_0))\cdot L_Y=-2+\frac{1}{2}\lambda+\sum_{i=1}^u\frac{\mu_i-1}{\mu_i}+\sum_{i=1}^n\sum_{j=1}^m\frac{w_{i,j}(C_i\cdot L)}{2}s_j(t_0).$$
We consider the g-pair
$$\left(\mathbb P^1,B_{\mathbb P^1}(t):=\sum_{i=1}^{u}\frac{\mu_i-1}{\mu_i}P_{i},\Mm(t):=\frac{\lambda}{2}\bar Q_0+\sum_{i=1}^n\sum_{j=1}^m\frac{w_{i,j}(C_i\cdot L)}{2}s_j(t)\bar Q_{i,j}\right)$$
where $P_i,Q_0,Q_{i,j}$ are distinct points on $\mathbb P^1$. Then $(\mathbb P^1,B_{\mathbb P^1}(t_0),\Mm(t_0))$ is lc, $K_{\mathbb P^1}+B_{\mathbb P^1}(t_0)+\Mm(t_0)_{\mathbb P^1}\equiv0$, and $K_{\mathbb P^1}+B_{\mathbb P^1}(r_c)+\Mm(r_c)_{\mathbb P^1}\not\equiv0$. The proposition now follows from \cite[Theorem 3.6]{Che23a} (which is essentially \cite[Theorem 3.8, Corollary 3.9]{Nak16} but the latter does not discuss the dimension $1$ case).
\end{proof}

\subsection{Threefold rank two case}\label{subsec: 3fold rk2}

In this subsection, we prove a slightly weaker version of Theorem \ref{thm: uniform rational polytope foliation intro} when $\dim X=3$ and $\rk\Ff=2$.

\begin{prop}\label{prop: uniform rational polytope foliation one variable threefold cork1}
Let $c,m$ be positive integers, $r_1,\dots,r_c$ real numbers such that $1,r_1,\dots,r_c$ are linearly independent over $\mathbb Q$, $\bm{r}:=(r_1,\dots,r_c)$, and $s_1,\dots,s_m: \mathbb R^{c+1}\rightarrow\mathbb R$ $\mathbb Q$-linear functions. Then there exists a positive real number $\delta$ depending only on $\bm{r}$ and $s_1,\dots,s_m$ satisfying the following. Assume that
\begin{enumerate}
    \item $(X,\Ff,B=\sum_{i=1}^ms_i(1,\bm{r})B_i)$ is an lc foliated triple such that $\dim X=3$ and $\rk\Ff=2$,
    \item $B_i\geq 0$ are distinct Weil divisors (possibly $0$) and $s_i(1,\bm{r})\geq 0$, and
    \item $B(t):=\sum_{i=1}^ms_i(1,r_1,\dots,r_{c-1},t)B_i$ for any $t\in\mathbb R$.
\end{enumerate}
Then $(X,\Ff,B(t))$ is lc for any $t\in (r_c-\delta,r_c+\delta)$.
\end{prop}
\begin{proof}
\noindent\textbf{Step 1}. In this step we apply Lemma \ref{lem: find nontrivial divisor on dlt model} and reduce to the case when $(X,\Ff,B)$ is $\Qq$-factorial dlt with additional good properties. This step is very similar to the beginning of the proof of Proposition \ref{prop: uniform rational polytope foliation one variable threefold rk1}.

We let $s_i(t):=s_i(1,r_1,\dots,r_{c-1},t)$ for any $t\in\mathbb R$. If $s_i(r_c)=0$, then $s_i(t)=0$ for any $i$, so we may assume that $s_i(r_c)\not=0$ for any $i$. By Lemma \ref{lem: dlt not influence perturbation} and \cite[Theorem 5.7]{LLM23}, possibly replacing $(X,\Ff,B(r_c))$ with a dlt model, we may assume that $(X,\Ff,B(r_c))$ is $\Qq$-factorial lc and $(X,\Ff,0)$ is dlt. In particular, $\Ff$ is non-dicritical. (We remark that $(X,\Ff,B(r_c))$ is actually dlt here, but later we will replace the foliated triple $(X,\Ff,B(r_c))$ again and it may no longer be dlt).

We only need to prove that there exists a positive real number $\epsilon$ depending only on $\Ii$ and $r$, such that for any lc threshold $t_0$ of $(X,\Ff,B(t))$, $|t_0-r_c|>\epsilon$. Thus we may assume that $(X,\Ff,B(t))$ has an lc threshold $t_0$. Since $1,r_1,\dots,r_c$ are linearly independent over $\Qq$, $r_c\not=t_0$. Moreover, there exists a positive real number $\delta_1$, such that for any $t\in (r_c-\delta_1,r_c+\delta_1)$ and any $i$, $s_i(t)\geq\frac{1}{2}s_i(r_c)>0$. In particular, for any $t\in (r_c-\delta_1,r_c+\delta_1)$, $\Supp B(t)=\Supp B(r_c)$, $\Supp\lfloor B(t)\rfloor=\Supp\lfloor B(r_c)\rfloor$, and $B(t)\geq\frac{1}{2}B(r_c)$.

We may assume that $t_0\in (r_c-\delta_1,r_c+\delta_1)$. In particular, we may assume that $(X,\Ff,B(t_0))$ does not have an lc center in codimension $1$ that is also an lc center of $(X,\Ff,B(r_c))$. Thus $(X,\Ff,B(t_0))$ has an lc center $x$ such that $\dim\bar x\leq 1$, and $x$ is not an lc center of $(X,\Ff,B(r_c))$. In particular, $x\in\Supp B(r_c)$.

By Lemma \ref{lem: find nontrivial divisor on dlt model} (more precisely, $X,\Ff,B,M$ in Lemma \ref{lem: find nontrivial divisor on dlt model} correspond to our $X,\Ff,B(r_c),B(t_0)-B(r_c)$ respectively), possibly replacing $(X,\Ff,B(t))$, we may assume that there exist a divisorial contraction $g: Y\rightarrow X$ of a prime divisor $E$ and a real number $s$ satisfying the following: let $B_Y(t)$ be the strict transform of $B(t)$ on $Y$ for any $t$ and $\Ff_Y:=g^{-1}\Ff$, then
\begin{itemize}
\item $s\in (r_c,t_0)$ if $r_c>t_0$, and $s\in (t_0,r_c)$ if $t_0<r_c$,
\item $(X,\Ff,B(s))$ is $\Qq$-factorial dlt, $(X,\Ff,B(r_c))$ is lc, and  $(X,\Ff,B(t_0))$ is lc,
\item $-E$ is ample over $X$,
\item $(Y,\Ff_Y,B_Y(s)+\epsilon_{\Ff}(E)E)$ is $\Qq$-factorial dlt, and
\item $a(E,\Ff,B(t_0))=-\epsilon_{\Ff}(E)$ and $a(E,\Ff,B(r_c))>-\epsilon_{\Ff}(E)$. In particular,  $(Y,\Ff_Y,B_Y(t_0)+\epsilon_{\Ff}(E)E)$ is lc.
\end{itemize}

\medskip

\noindent\textbf{Step 2}. We deal with the case when $\epsilon_{\Ff}(E)=1$.

Suppose that $\epsilon_{\Ff}(E)=1$. Since $(Y,\Ff_Y,B_Y(s)+\epsilon_{\Ff}(E))$ is $\Qq$-factorial dlt, $\Ff_Y$ is non-dicritical. By \cite[Remark 2.16]{CS21} and \cite[Lemma 3.11]{Spi20}, over the generic point of $\Center_XE$,
\begin{itemize}
    \item $(X,B(r_c))$ is lc, and $(X,B(t_0))$ is lc, and
    \item $a(E,X,B(t_0))=-1$ and $a(E,X,B(r_c))>-1$.
\end{itemize}
The proposition follows from \cite[Theorem 1.6]{Nak16} (see also \cite[Theorem 5.6]{HLS19}) in this case.

\medskip

\noindent\textbf{Step 3}. From now on we can assume that $\epsilon_{\Ff}(E)=0$. We summarize some known properties in this step.

Let $V$ be the normalization of $\Center_XE$, $E^\nu$ the normalization of $E$,  $g|_{E^\nu}: E^\nu\rightarrow V$ the induced projective surjective morphism with Stein factorization
$$E^\nu\xrightarrow{\pi} W\xrightarrow{\tau} V,$$
and $F$ a general fiber of $\pi$. Then since $-E$ is ample over $X$, 
$$K_{\Ff_Y}+B_Y(r_c)=g^*(K_{\Ff}+B(r_c))+a(E,\Ff,B(r_c))E\sim_{\mathbb R,X}a(E,\Ff,B(r_c))E$$
is anti-ample$/X$. Thus
$$(K_{\Ff_Y}+B_Y(r_c))|_F$$
is anti-ample. Moreover, since
$$K_{\Ff_Y}+B_Y(t_0)=g^*(K_{\Ff}+B(t_0))\sim_{\mathbb R,X}0,$$
we have that
$$(K_{\Ff_Y}+B_Y(t_0))|_F\sim_{\mathbb R}0.$$
By Theorem \ref{thm: adjunction formula corank 1 dim 3 intro}, there exist positive integers $w_1,\dots,w_n$, $\lambda_1,\dots,\lambda_k$, non-negative integers $\{w_{i,j}\}_{1\leq i\leq n, 1\leq j\leq m}$, and prime divisors $C_1,\dots,C_n,T_1,\dots,T_l$, such that
for any real number $t$, we have
$$\left(K_{\Ff_Y}+B_Y(t)\right)|_{E^\nu}=K_{E^\nu}+\sum_{i=1}^n\frac{w_i-1+\sum_{j=1}^m w_{i,j}s_j(t)}{w_i}C_i+\sum_{i=1}^l\lambda_iT_i.$$
We let $$B_E(t):=\sum_{i=1}^n\frac{w_i-1+\sum_{j=1}^m w_{i,j}s_j(t)}{w_i}C_i+\sum_{i=1}^l\lambda_iT_i$$ and
$$\tilde B_E(t):=\sum_{i=1}^n\frac{w_i-1+\sum_{j=1}^m w_{i,j}s_j(t)}{w_i}C_i+\sum_{i=1}^lT_i$$ for any real number $t$. By Theorem \ref{thm: adjunction formula corank 1 dim 3 intro}(3), $(E^\nu,\tilde B_E(t_0))$ and  $(E^\nu,\tilde B_E(r_c))$ are lc.

\medskip

\noindent\textbf{Step 4}. We deal with the case when $\dim F=1$. 

Suppose that $\dim F=1$. We let 
$$h(t):=B_E(t)\cdot F=\sum_{i=1}^n\frac{w_i-1+\sum_{j=1}^m w_{i,j}s_j(t)}{w_i}(C_i\cdot F)+\sum_{i=1}^l\lambda_i(T_i\cdot F)$$
for any real number $t$. Then $h(t)$ is an affine function, $h(r_c)\geq 0$, and $h(t_0)\geq 0$. 

In this case, we have $K_{E^\nu}\cdot F+h(t_0)=0$ and $K_{E^\nu}\cdot F+h(r_c)<0$. Thus $K_{E^\nu}\cdot F<0$, so $F$ is a smooth rational curve. Since $F^2=0$, $K_{E^\nu}\cdot F=-2$. Thus
$$-2+\sum_{i=1}^n\frac{w_i-1+\sum_{j=1}^m w_{i,j}s_j(t)}{w_i}(C_i\cdot F)+\sum_{i=1}^l\lambda_i(T_i\cdot F)=0.$$
We consider the pair $$\left(\mathbb P^1,B_{\mathbb P^1}(t):=\sum_{i=1}^n\frac{w_i-1+\sum_{j=1}^mw_{i,j}s_j(t)}{w_i}\sum_{k=1}^{(C_i\cdot F)}P_{i,k}+\sum_{i=1}^l\sum_{k=1}^{\lambda_i(T_i\cdot F)}Q_{i,k}\right)$$
for any real number $t$, where $P_{i,k}$ are $Q_{i,k}$ are distinct points on $\mathbb P^1$. Then $(\mathbb P^1,B_{\mathbb P^1}(t_0))$ is lc, $K_{\mathbb P^1}+B_{\mathbb P^1}(t_0)\equiv0$, and $K_{\mathbb P^1}+B_{\mathbb P^1}(r_c)\not\equiv 0$, The proposition now follows from \cite[Theorem 3.6]{Che23a}.

\medskip

\noindent\textbf{Step 5}. We conclude the proof in this step. Now we may assume that $F=E^\nu$ and $\dim F=2$. Then $K_{E^\nu}+B_E(t_0)\sim_{\mathbb R}0$, and $K_{E^\nu}+B_E(r_c)$ is anti-ample. By  \cite[Theorem 3.8, Corollary 3.9]{Nak16}, we may assume that $B_E(t)\not=\tilde B_E(t)$ for any $t$. Thus $c_i\geq 2$ for some $i$. We let $T_i^\nu$ be the normalization of $T_i$.

By Theorem \ref{thm: adjunction formula corank 1 dim 3 intro}, there exist a non-negative integer $n_i$, positive integers $e_{i,1},\dots,e_{i,n_i}$, non-negative integers $\{e_{i,k,j}\}_{1\leq k\leq n_i,1\leq j\leq m}$, a Weil divisor $P_i\geq 0$ on $T_i^\nu$, and prime divisors $P_{i,1},\dots,P_{i,n_i}$ on $T_i^\nu$, such that for any real number $t$,
$$(K_{\Ff_Y}+B_Y(t))|_{T_i^\nu}=K_{T_i^\nu}+B_{T_i}(t):=K_{T_i^\nu}+P_i+\sum_{k=1}^{n_i}\frac{e_{i,k}-1+\sum_{j=1}^me_{i,k,j}s_j(t)}{e_{i,k}}P_{i,k},$$
and $$\left(T_i^\nu,\Supp P_i+\sum_{k=1}^{n_i}\frac{e_{i,k}-1+\sum_{j=1}^me_{i,k,j}s_j(t_0)}{e_{i,k}}P_{i,k}\right)$$ is lc. Moreover, we have $K_{T_i^\nu}+B_{T_i}(t_0)\sim_{\mathbb R}0$ and $K_{T_i^\nu}+B_{T_i}(r_c)$ is anti-ample. Thus $T_i^\nu=\mathbb P^1$. We consider the g-pair
$$\left(\mathbb P^1,B'_{\mathbb P^1}(t):=\sum_{k=1}^{n_i}\frac{e_{i,k}-1+\sum_{j=1}^me_{i,k,j}s_j(t)}{e_{i,k}}P_{i,k},\Mm:=\bar P_i\right)$$
 for any real number $t$. Then $(\mathbb P^1,B'_{\mathbb P^1}(t_0),\Mm)$ is lc, $K_{\mathbb P^1}+B'_{\mathbb P^1}(t_0)+\Mm_{\mathbb P^1}\equiv 0$, and $K_{\mathbb P^1}+B'_{\mathbb P^1}(r_c)+\Mm_{\mathbb P^1}\not\equiv 0$. The proposition now follows from \cite[Theorem 3.6]{Che23a}.
\end{proof}

\begin{rem}
The same arguments as of Proposition \ref{prop: uniform rational polytope foliation one variable threefold rk1} can also be helpful to simplify the proof of the rank $2$ case of the ACC for foliated lc thresholds for threefolds \cite[Theorem 3.11]{Che22}. By applying Lemma \ref{lem: find nontrivial divisor on dlt model}, we do not need the argument from  \cite[Page 15]{Che22} to \cite[End of Section 3]{Che22}.
\end{rem}

\subsection{Proof of Theorem \ref{thm: uniform rational polytope foliation intro}}

\begin{thm}\label{thm: uniform rational polytope foliation one variable}
Let $c,m$ be positive integers, $r_1,\dots,r_c$ real numbers such that $1,r_1,\dots,r_c$ are linearly independent over $\mathbb Q$, $\bm{r}:=(r_1,\dots,r_c)$, and $s_1,\dots,s_m: \mathbb R^{c+1}\rightarrow\mathbb R$ $\mathbb Q$-linear functions. Then there exists a positive real number $\delta$ depending only on $\bm{r}$ and $s_1,\dots,s_m$ satisfying the following. Assume that
\begin{enumerate}
    \item $(X,\Ff,B=\sum_{i=1}^ms_i(1,\bm{r})B_i)$ is an lc foliated triple such that $\rk\Ff<\dim X\leq 3$,
    \item $B_i\geq 0$ are distinct Weil divisors (possibly $0$) and $s_i(1,\bm{r})\geq 0$,
    \item $B(t):=\sum_{i=1}^ms_i(1,r_1,\dots,r_{c-1},t)B_i$ for any $t\in\mathbb R$, and
    \item if $\dim X=3$ and $\rk\Ff=1$, then $X$ is projective.
\end{enumerate}
Then $(X,\Ff,B(t))$ is lc for any $t\in (r_c-\delta,r_c+\delta)$.
\end{thm}
\begin{proof}
If $\rk\Ff=0$ then $B(t)=0$ for any $t$ and there is nothing left to prove, so we may assume that $\rk\Ff>0$. The theorem follows from Propositions \ref{prop: uniform rational polytope foliation one variable threefold rk1} and \ref{prop: uniform rational polytope foliation one variable threefold cork1} and \cite[Theorem 1.8]{LMX24a}.
\end{proof}

\begin{proof}[Proof of Theorem \ref{thm: uniform rational polytope foliation intro}]
If $\Ff=T_X$, then the theorem is \cite[Theorem 5.6]{HLS19}. So we may assume that $\rk\Ff<\dim X$.

We apply induction on $c$. When $c=1$, Theorem \ref{thm: uniform rational polytope foliation intro} directly follows from Theorem \ref{thm: uniform rational polytope foliation one variable}. When $c\geq 2$, by Theorem \ref{thm: uniform rational polytope foliation one variable}, there exists a positive integer $\delta$ depending only on $r_1,\dots,r_c,s_1,\dots,s_m$, such that for any $t\in (r_c-\delta,r_c+\delta)$, $(X,\Ff,\sum_{i=1}^ms_i(1,r_1,\dots,r_{c-1},t)B_i)$ is lc. We pick rational numbers $r_{c,1}\in (r_c-\delta,r_c)$ and $r_{c,2}\in (r_c,r_c+\delta)$ depending only on $r_1,\dots,r_c,s_1,\dots,s_m$. By induction on $c$, there exists an open subset $U_0\ni (r_1,\dots,r_{c-1})$ of $\mathbb R^{c-1}$, such that for any $\bm{v}\in U_0$, $(X,\Ff,\sum_{i=1}^ms_i(1,\bm{v},r_{c,1})B_i)$ and $(X,\Ff,\sum_{i=1}^ms_i(1,\bm{v},r_{c,2})B_i)$ are lc. We may pick $U:=U_0\times (r_{c,1},r_{c,2})$.
\end{proof}

\section{Proofs of Theorem \ref{thm: global acc threefold} and Corollary \ref{cor: accumulation point foliated threefold lct}}\label{sec: proof of the main theorems}

The following result is a variation of the ACC for lc thresholds \cite[Theorem 0.5]{Che22} for foliations in dimension $\leq 3$, which is more useful in some scenarios.

\begin{prop}\label{prop: variation acc lct foliation}
Let $\Ii\subset [0,1]$ be a DCC set of real numbers and $\alpha$ a positive real number. Then there exists a function $g: \bar\Ii\rightarrow\bar\Ii$ satisfying the following.
\begin{enumerate}
\item $g\circ g=g$ and $g(\bar\Ii)$ is a finite set.
    \item $\gamma+\alpha\geq g(\gamma)\geq \gamma$ for any $\gamma\in\bar\Ii$.
    \item $g(\gamma)\leq g(\gamma')$ for any $\gamma,\gamma'\in\bar\Ii$ such that $\gamma\leq\gamma'$.
    \item For any non-negative integer $m$ and lc foliated triple $(X,\Ff,\sum_{i=1}^mb_iB_i)$, such that $\dim X\leq 3$, $b_i\in\Ii$ for any $i$, and $B_i$ are effective $\Qq$-Cartier Weil divisors, we have that $$\left(X,\Ff,\sum_{i=1}^mg(b_i)B_i\right)$$
    is lc.
\end{enumerate}
\end{prop}
\begin{proof}
We may assume that $\Ii=\bar\Ii$. Let
$$\Ii':=\overline{\{\lct(X,\Ff,B;D)\mid \dim X\leq 3, X\text{ is projective}, B\in\bar\Ii,D\in\mathbb N^+\}}.$$
By \cite[Theorem 0.5]{Che22}, $\Ii'$ is an ACC set. By \cite[Lemma 5.17]{HLS19}, there exists a function $g:\bar\Ii\rightarrow\bar\Ii$, such that (1-3) hold, and for any $\beta\in\Ii'$ and $\gamma\in\bar\Ii$ such that $\beta\geq\gamma$, we have $\beta\geq g(\gamma)$.

Suppose that $(X,\Ff,\sum_{i=1}^mg(b_i)B_i)$ is not lc. Then $m\geq 1$, and there exists $0\leq j\leq m-1$, such that $(X,\Ff,\sum_{i=1}^jg(b_i)B_i+\sum_{i=j+1}^mb_iB_i)$ is lc and  $(X,\Ff,\sum_{i=1}^{j+1}g(b_i)B_i+\sum_{i=j+2}^mb_iB_i)$ is not lc. Let
$$b:=\lct\left(X,\Ff,\sum_{i\not=j+1}g(b_i)B_i;B_{j+1}\right),$$
then $b_{j+1}\leq b<g(b_{j+1})$ and $b\in\Ii'$, which is not possible. Thus $(X,\Ff,\sum_{i=1}^mg(b_i)B_i)$ is lc and we are done.
\end{proof}

\begin{prop}\label{prop: global acc with contraction}
Let $\Ii\subset [0,1]$ be a DCC set. Then there exist a finite set $\Ii_0\subset\Ii$ depending only on $\Ii$ satisfying the following. Assume that
\begin{enumerate}
    \item $(X,\Ff,B)$ is a projective lc foliated triple of dimension $\leq 3$,
    \item $B\in\Ii$,
    \item $\pi: X\rightarrow Z$ is a contraction such that $0<\dim Z<\dim X$,
    \item  there exists a foliation $\Ff_Z$ on $Z$ such that $\Ff=\pi^{-1}\Ff_Z$, 
    \item $K_{\Ff}+B\sim_{\mathbb R}0$, and
    \item if $\dim X=3,\rk\Ff=2$, and $\dim Z=2$, then $Z$ is klt.
\end{enumerate} 
Then $B\in\Ii_0$.
\end{prop}
In the proof of Proposition \ref{prop: global acc with contraction}, we will introduce a lot of sets of coefficients depending only on $\Ii$. For the reader's convenience, we note that the sets with/without the subscript ``$0$" are finite/DCC sets in the following proof.

\begin{proof}
By \cite[Theorem 1.5]{HMX14}, we may assume that $\rk\Ff<\dim X$. If $\rk\Ff=0$, then $B=0$ and there is nothing left to prove. So we may assume that $1\leq\rk\Ff\leq 2$. 

By \cite[Theorem 5.7]{LLM23}, possibly replacing $\Ii$ with $\Ii\cup\{1\}$, we may assume that $(X,\Ff,B)$ is a $\Qq$-factorial projective dlt foliated triple of dimension $\leq 3$.

 Let $F$ be a general fiber of $\pi$ and $B_F:=B|_F$. Then
 $$K_F+B_F=(K_X+B)|_F=(K_{\Ff}+B)|_F\sim_{\mathbb R}0,$$
 and $B_F\in\Ii$. Since $(X,\Ff,B)$ is lc, $(F,B_F)$ is lc. By \cite[Theorem 1.5]{HMX14}, there exists a finite set $\Ii_0'\subset\Ii$ such that $B_F\in\Ii_0'$.
 
 By \cite[Theorem 2.5]{Che22}, we may assume that $\dim X=3$. By \cite[Proposition 6.4(2)]{LLM23}, there exists a projective lc generalized foliated quadruple $(Z,\Ff_Z,B_Z,\Mm)$ induced by a canonical bundle formula $\pi: (X,\Ff,B)\rightarrow Z$ (see \cite[Definition 1.2]{LLM23} for the definition) Then $K_{\Ff_Z}+B_Z+\Mm_Z\sim_{\mathbb R}0$. 

We let $B^h$ be the horizontal$/Z$ part of $B$ and $B^v$ the vertical$/Z$ part of $B$. Then $B^h\in\Ii_0'$. Let $\bar \Ii$ be the closure of $\Ii$. By Proposition \ref{prop: variation acc lct foliation}, there exist a finite set $\Ii_0''\subset\bar\Ii$ depending only on $\Ii$ and an $\Rr$-divisor $\bar B\in\Ii_0''$, such that \begin{itemize}
\item $\Ii_0'\subset\Ii_0''$,
    \item $\bar B\geq B$ and $\Supp\bar B=\Supp B$,
    \item $(X,\Ff,\bar B)$ is lc, and
    \item $\bar B^h=B^h$, where $\bar B^h$ is the horizontal$/Z$ part of $\bar B$.
\end{itemize}
In particular, $0\leq \bar B-B$ is vertical$/Z$. 

\begin{claim}\label{claim: prop 8.1 rk2 algint}
Proposition \ref{prop: global acc with contraction} holds if  $\dim Z=1$.
\end{claim}
\begin{proof}
In this case, $\rk\Ff=2$ and $\Ff_Z$ is the trivial foliation, so $K_{\Ff_Z}+B_Z+\Mm_Z\sim_{\mathbb R}0$ implies that $B_Z=0$. Suppose that there exists a component $D$ of $B$ that is vertical over $Z$. We let $P=\pi(D)$, then 
$$\sup\{t\mid (X,\Ff,B+t\pi^*P)\text{ is sub-lc over }P\}<1.$$
By \cite[Proposition 6.4(2)]{LLM23}, $\mult_PB_Z>0$, a contradiction. Therefore, all components of $B$ are horizontal over $Z$. Since $B^h\in\Ii_0'$, Proposition \ref{prop: global acc with contraction} follows in this case.
\end{proof}

\noindent\textit{Proof of Proposition \ref{prop: global acc with contraction} continued}. By Claim \ref{claim: prop 8.1 rk2 algint}, we may assume that $\dim Z=2$. Since $(X,\Ff,\bar B)$ is lc, $(X,\Ff,B)$ is dlt, and $\Supp\bar B=\Supp B$, we have that $(X,\Ff,\frac{1}{2}(\bar B+B))$ is dlt. By \cite[Theorem 5.8]{LLM23}, we may run a $(K_{\Ff}+\frac{1}{2}(\bar B+B))$-MMP$/Z$ which terminates with a log minimal model  (cf. \cite[Definition 5.5]{LLM23})  of $(X,\Ff,\frac{1}{2}(\bar B+B))/Z$.  Since
$$\frac{1}{2}(K_{\Ff}+\bar B)\sim_{\mathbb R}K_{\Ff}+\frac{1}{2}(\bar B+B),$$
this MMP is also a $(K_{\Ff}+\bar B)$-MMP which terminates with a weak lc model  (cf. \cite[Definition 5.5]{LLM23}) $(X',\Ff',\bar B')/Z$ of $(X,\Ff,\bar B)/Z$. We let $\pi': X'\rightarrow Z$ be the induced contraction and $B'$  the image of $B$ on $X'$.

\begin{claim}\label{claim: kf'+b' is trivial over Z}
$K_{\Ff'}+\bar B'\sim_{\mathbb R,Z}0$.
\end{claim}
\begin{proof}
Since $K_{\Ff'}+B'\sim_{\mathbb R}0$, we only need to show that $\bar B'-B'\sim_{\mathbb R,Z}0$. By assumption, $\bar B'-B'$ is nef$/Z$. There are two cases.

\medskip

\noindent\textbf{Case 1}. $\rk\Ff=1$. In this case, $\Ff_Z=0$, so $K_{\Ff_Z}+B_Z+\Mm_Z\sim_{\mathbb R}0$ implies that $B_Z=0$. Therefore, for any component $D$ of $B$ that is vertical over $Z$, $D$ is very exceptional over $Z$. Thus $\bar B-B$ is very exceptional over $Z$, so $\bar B'-B'$ is very exceptional over $Z$. By \cite[Lemma 3.3]{Bir12}, $\bar B-B=0$ and we are done.

\medskip

\noindent\textbf{Case 2}. $\rk\Ff=2$. In this case, by our assumption, $Z$ is klt. Thus $Z$ is $\Qq$-factorial. For any prime divisor $D$ on $Z$, we define 
$$\nu_D:=\sup\{t\mid \bar B'-B'-\pi^*D\geq 0\}$$
and let 
$$L':=\bar B'-B'-\sum_{D\text{ is a prime divisor on }Z}\nu_D\pi'^*D.$$
Then $L'\geq 0$ and $L'$ is very exceptional over $Z$.  By \cite[Lemma 3.3]{Bir12}, $L'=0$. Thus $$\bar B'-B'=\sum_{D\text{ is a prime divisor on }Z}\nu_D\pi'^*D\sim_{\mathbb R,Z}0$$ 
and we are done.
\end{proof}

\noindent\textit{Proof of Proposition \ref{prop: global acc with contraction} continued}. By Claim \ref{claim: kf'+b' is trivial over Z}, $K_{\Ff'}+\bar B'\sim_{\mathbb R,Z}0$. Since $K_{\Ff'}+B'\sim_{\mathbb R}0$, $\bar B'-B'\sim_{\mathbb R,Z}0$. By Theorem \ref{thm: uniform rational polytope foliation intro} and \cite[Lemma 5.3]{HLS19}, there exist real numbers $a_1,\dots,a_k\in (0,1]$, a finite set $\Ii_0'''\subset\mathbb Q\cap [0,1]$ depending only on $\Ii$, and $\Qq$-divisors $\bar B_1',\dots,\bar B_k'\in\Ii_0'''$ on $X'$, such that
\begin{itemize}
    \item $\sum_{i=1}^ka_i=1$ and $\sum_{i=1}^ka_i\bar B_i'=\bar B'$,
    \item $(X',\Ff',\bar B'_i)$ is lc for each $i$,
    \item $(X',\Ff',B'_i:=\bar B'_i-(\bar B'-B'))$ is lc for each $i$, and
    \item $K_{\Ff'}+\bar B'_i\sim_{\mathbb Q,Z}0$ and  $K_{\Ff'}+B'_i\sim_{\mathbb R,Z}0$ for each $i$.
\end{itemize}
By \cite[Proposition 6.4(2.a)(3)]{LLM23}, we may let $(Z,\Ff_{Z},B_{i,Z},\Mm_i)$ and $(Z,\Ff_{Z},\bar B_{i,Z},\Mm_i)$ be projective lc generalized foliated quadruples induces by canonical bundle formulas of $(X',\Ff',B'_i)\rightarrow Z$ and $(X',\Ff',\bar B'_i)\rightarrow Z$ respectively.
We let $\Mm':=\sum_{i=1}^ka_i\Mm_i$ (note that it is possible that $\Mm'\not=\Mm$), and let $B'_{Z}:=\sum_{i=1}^ka_iB_{i,Z}$. Then $K_{\Ff_{Z}}+B'_{Z}+\Mm'_{Z}\sim_{\mathbb R}0$. Since $\dim Z=2$, by \cite[Proposition 6.4(4)]{LLM23}, there exists a positive integer $I$ depending only on $\Ii$, such that for each $i$, we may choose $\Mm_i$ such that $I\Mm_i$ is base-point-free, $I(K_{\Ff'}+B_i')\sim I\pi'^*(K_{\Ff_{Z}}+B_{i,Z}+\Mm_{i,Z})$, and $I(K_{\Ff'}+\bar B_i')\sim I\pi'^*(K_{\Ff_{Z}}+\bar B_{i,Z}+\Mm_{i,Z})$.

By \cite[Theorem 0.5]{Che22} and \cite[Proposition 6.4(2)]{LLM23}, there exists a DCC set $\Ii'\subset [0,1]$ depending only on $\Ii$ such that $\bar B_{i,Z},B_{i,Z}\in\Ii'$. Since
$K_{\Ff_{Z}}+B_{i,Z}+\Mm_{i,Z}\sim_{\mathbb R}0$,
by \cite[Lemma 7.2]{LLM23}, there exists a finite set $\Ii_0''''$ depending only on $\Ii$ such that $B'_{i,Z}\in\Ii_0''''$ for any $i$. Possibly replacing $\Ii_0''''$, we may assume that $B_Z'\in\Ii_0''''$.

We let $\bar B_{Z}:=\sum_{i=1}^ka_i\bar B_{i,Z}$. Then there exists a DCC set $\Ii''\subset (0,1]$ depending only on $\Ii$ such that for any component $D$ of $\bar B_{Z}-B'_{Z}$, $\mult_D(\bar B_{Z}-B'_{Z})\in\Ii''$. Since the coefficients of $\bar B'$ belongs to the finite set $\Ii_0''$, the coefficients of
$$B'=\bar B'-\pi'^*(\bar B_{Z}-B'_{Z})$$
belong to an ACC set depending only on $\Ii$. Since $B'\in\Ii$, the coefficients of $B'$ belong to a finite set $\tilde\Ii_0$ depending only on $\Ii$.  In this case, by Theorem \ref{thm: uniform rational polytope foliation intro} and \cite[Lemma 5.3]{HLS19}, there exist real numbers $c_1,\dots,c_k\in (0,1]$ and a finite set $\tilde\Ii_0'\subset\mathbb Q\cap [0,1]$ depending only on $\Ii$, and $\Qq$-divisors $\tilde B_1',\dots,\tilde B_k'\in\tilde\Ii_0'$ on $X'$, such that
\begin{itemize}
    \item $\sum_{i=1}^kc_i=1$ and $\sum_{i=1}^kc_i\tilde B_i'=B'$,
    \item $(X',\Ff',\tilde B'_i)$ is lc for each $i$, and
    \item $K_{\Ff'}+\tilde B'_i\sim_{\mathbb Q}0$ for each $i$.
\end{itemize}
By \cite[Theorem 1.3]{LMX24a}, there exists a positive integer $I$ depending only on $\Ii$ such that $I(K_{\Ff'}+\tilde B_i')\sim 0$ for each $i$. Thus for any prime divisor $E$ over $X'$, $a(E,\Ff',\tilde B_i')$ belong to the discrete set $$\left\{\frac{k}{I}\biggm| k\geq -I, k\in\mathbb Z\right\}.$$ In particular, for any component $D$ of $B$,
$$\mult_{D}B=-a(D,\Ff,B)=-a(D,\Ff',B')=-\sum_{i=1}^kc_ia(D,\Ff',\tilde B_i')$$
belongs to a discrete set. Since $\mult_DB\in [0,1]$, $\mult_DB$ belongs to a finite set $\Ii_0$. The theorem follows.
\end{proof}

\begin{proof}[Proof of Theorem \ref{thm: global acc threefold}]
By \cite[Theorem 2.5]{Che22}, we may assume that $\dim X=3$. By \cite[Theorem 5.7]{LLM23}, possibly replacing $\Ii$ with $\Ii\cup\{1\}$ and $(X,\Ff,B)$ with its dlt model, we may assume that $(X,\Ff,B)$ is $\Qq$-factorial dlt. If $B=0$, then there is nothing left to prove, so we may assume that $B\not=0$, hence $K_{\Ff}$ is not pseudo-effective. By \cite[Theorem 1.1]{CP19},\cite[Theorem 3.1]{LLM23}, there exists an algebraically integrable foliation $0\not=\Ee\subset\Ff$. If $\Ff$ is algebraically integrable, then by \cite[Theorem 3.10]{ACSS21}, possibly replacing $(X,\Ff,B)$, we may assume that there exists a contraction $\pi: X\rightarrow Z$ such that $\Ff$ is induced by $\pi$ and $Z$ is smooth, and the theorem follows from Proposition \ref{prop: global acc with contraction}. Thus we may assume that $\Ff$ is not algebraically integrable. 

Let $S$ be a component of $B$ and let $b:=\mult_SB$. By \cite[Lemma 8.2]{LLM23}, there exists a birational map $f: X\dashrightarrow X'$ and a contraction $\pi': X'\rightarrow Z$ such that $f$ does not extract any divisor, $S$ is not contracted by $f$, $\Ff':=f_*\Ff$ is induced by a foliation $\Ff_Z$ on $Z$ (i.e., $\Ff'=\pi'^{-1}\Ff_Z$), $\dim Z=2$, and $Z$ is klt. Let $B':=f_*B$ and $S':=f_*S$. Then $S'\not=0$. Since $K_{\Ff}+B\sim_{\mathbb R}0$, $K_{\Ff'}+B'\sim_{\mathbb R}0$ and $(X',\Ff',B')$ is lc. By Proposition \ref{prop: global acc with contraction}, there exists a finite set $\Ii_0$ depending only on $\Ii$ such that $b\in\Ii_0$. Since $S$ can be any component of $B$, $B\in\Ii_0$, and we are done.
\end{proof}

\begin{thm}\label{thm: accumulation point of foliated lc threshold complicated version}
Let $c$ be a non-negative integer, $r_1,\dots,r_c$ real numbers, and $\Ii\subset [0,1]$ a DCC set, such that $\bar\Ii\subset\Span_{\mathbb Q}(1,r_1,r_2,\dots,r_c)$. The the accumulation points of 
$$\{\lct(X,\Ff,B;D)\mid \dim X\leq 3, (X,\Ff,B)\text{ is lc, }B\in\Ii,D\in\mathbb N^+\}$$
belong to $\Span_{\mathbb Q}(1,r_1,r_2,\dots,r_c)$.
\end{thm}
\begin{proof}
Suppose the theorem does not hold. By \cite[Theorem 0.5]{Che22}, there exist a sequence of lc foliated triples $(X_i,\Ff_i,B_i)$ of dimension $\leq 3$ and effective $\Qq$-Cartier Weil divisors $D_i$ on $X_i$, such that $t_i:=\lct(X_i,\Ff_i,B_i;D_i)$ is strictly decreasing, $B_i\in\Ii$, and $t:=\lim_{i\rightarrow+\infty}t_i\not\in\Span_{\mathbb Q}(1,r_1,r_2,\dots,r_c)$. We write $B_i=\sum_{j=1}^{m_i}b_{i,j}B_{i,j}$, where $B_{i,j}$ are the irreducible components of $B_i$. Possibly replacing $(X_i,\Ff_i,B_i)$ with a dlt model and replacing $D_i$ with its pullback, we may assume that $(X_i,\Ff_i,B_i)$ is $\Qq$-factorial dlt. 

Since $(X_i,\Ff_i,B_i+t_iD_i)$ is lc, $(X_i,\Ff_i,B_i+tD_i)$ is lc, and the coefficients of $B_i+tD_i$ belong to a DCC set depending only on $\Ii$. By Proposition \ref{prop: variation acc lct foliation}, there exists a function $g: \bar\Ii\rightarrow\bar\Ii$, such that
\begin{enumerate}
\item $g\circ g=g$ and $\Ii_0:=g(\bar\Ii)$ is a finite set,
\item $g(\gamma)\geq \gamma$ for any $\gamma\in\bar\Ii$,
\item $g(\gamma)\leq g(\gamma')$ for any $\gamma,\gamma'\in\bar\Ii$ such that $\gamma\leq\gamma'$, and
\item  $$\left(X_i,\Ff_i,\sum_{j=1}^{m_i}g(b_{i,j})B_{i,j}+tD_i\right)$$
is lc for any $i$.
\end{enumerate}
Since $\Ii_0$ is a finite set, we have $\Ii_0=\{\bar b_1,\dots,\bar b_m\}$ for some non-negative integer $m$, and we may write $\sum_{j=1}^{m_i}g(b_{i,j})B_{i,j}=\sum_{j=1}^m\bar b_jC_{i,j}$ where $C_{i,j}$ are effective Weil divisors.

By our assumption, $\bar b_j\in\Span_{\mathbb Q}(1,r_1,\dots,r_c)$ for each $j$ and $t\not\in\Span_{\mathbb Q}(1,r_1,\dots,r_c)$. We let $V$ be the rational envelope of $(\bar b_1,\dots,\bar b_m,t)$ in $\mathbb R^{m+1}$, then $V=V'\times\mathbb R$, where $V'$ is the rational envelope of $(\bar b_1,\dots,\bar b_m)$ in $\mathbb R^m$. By Theorem \ref{thm: uniform rational polytope foliation intro}, there exist an open subset $U'\ni (\bar b_1,\dots,\bar b_m)$ of $V'$, and an open subset $W\ni t$ of $\mathbb R$, such that $(X_i,\Ff_i,\sum_{j=1}^mv_jC_{i,j}+wD_i)$ is lc for any $(v_1,\dots,v_m)\in U'$ and $w\in W$. In particular, $$\left(X_i,\Ff_i,\sum_{j=1}^m\bar b_jC_{i,j}+wD_i\right)=\left(X_i,\Ff_i,\sum_{j=1}^{m_i}g(b_{i,j})B_{i,j}+wD_i\right)$$ is lc for any $w\in W$. Possibly passing to a subsequence, we may assume that there exists a real number $w_0\in W$ such that $w_0>t_i$ for any $i$. Then $(X_i,\Ff_i,\sum_{j=1}^{m_i}g(b_{i,j})B_{i,j}+w_0D_i)$ is lc for any $i$, so $(X_i,\Ff_i,B_i+w_0D_i)$ is lc for any $i$, so
$$t_i=\lct(X_i,\Ff_i,B_i;D_i)\geq w_0>t_i,$$
a contradiction.
\end{proof}

\begin{proof}[Proof of Corollary \ref{cor: accumulation point foliated threefold lct}]
It immediately follows from Theorem \ref{thm: accumulation point of foliated lc threshold complicated version} by taking $c=0$.
\end{proof}


\begin{thebibliography}{99}

\bibitem[ACSS21]{ACSS21} F. Ambro, P. Cascini, V. V. Shokurov, and C. Spicer, \textit{Positivity of the moduli part}, arXiv:2111.00423.

\bibitem[AD14]{AD14} C. Araujo and S. Druel, \textit{On codimension 1 del Pezzo foliations on varieties with mild singularities}, Math. Ann., \textbf{360} (2014), no. 3--4, 769--798.


\bibitem[Bir12]{Bir12} C. Birkar, \textit{Existence of log canonical flips and a special LMMP}, Pub. Math. IHES., \textbf{115} (2012), 325--368.

\bibitem[Bru02]{Bru02} M. Brunella, \textit{Foliations on complex projective surfaces}, arXiv:math/0212082.

\bibitem[Bru15]{Bru15} M. Brunella, \textit{Birational geometry of foliations}, IMPA Monographs \textbf{1} (2015), Springer, Cham.

\bibitem[BCHM10]{BCHM10}
C. Birkar, P. Cascini, C. D. Hacon and J. M\textsuperscript{c}Kernan, \textit{Existence of minimal models for varieties of log general type}, J. Amer. Math. Soc. \textbf{23} (2010), no. 2, 405--468.

\bibitem[Che23a]{Che23a} G. Chen, \textit{Boundedness of n-complements for generalized pairs}, Eur. J. Math. \textbf{9} (2023), no. 95.

\bibitem[CHL23]{CHL23} G. Chen, J. Han, and J. Liu, \textit{On effective log Iitaka fibrations and existence of complements},  Int. Math. Res. Not. (2023), rnad253.

\bibitem[CHLX23]{CHLX23} G. Chen, J. Han, J. Liu, and L. Xie, \textit{Minimal model program for algebraically integrable foliations and generalized pairs}, arXiv:2309.15823.

\bibitem[Che22]{Che22} Y.-A. Chen, \textit{ACC for foliated log canonical thresholds}, arXiv:2202.11346.

\bibitem[Che23b]{Che23b} Y.-A. Chen, \textit{Log canonical foliation singularities on surfaces}, Math. Nachr. \textbf{296} (2023), no. 8, 3222--3256.


\bibitem[CP19]{CP19} F.~Campana and M.~P\u{a}un, \textit{Foliations with positive slopes and birational stability of orbifold cotangent bundles}, Publ. Math. Inst. Hautes \'{E}tudes Sci. \textbf{129} (2019), 1--49.

\bibitem[CS20]{CS20} P. Cascini and C. Spicer, \textit{On the MMP for rank one foliations on threefolds}, arXiv:2012.11433.

\bibitem[CS21]{CS21} P.~Cascini and C.~Spicer, \textit{MMP for co-rank one foliations on threefolds}, Invent. math. \textbf{225} (2021), 603--690.

\bibitem[CS23a]{CS23a} P. Cascini and C. Spicer, \textit{On the MMP for algebraically integrable foliations}, to appear in Shokurov's 70th birthday's special volume, arXiv:2303.07528.

\bibitem[CS23b]{CS23b} P. Cascini and C. Spicer, \textit{Foliation adjunction}, arXiv:2309.10697.

\bibitem[DLM23]{DLM23} O. Das, J. Liu, and R. Mascharak, \textit{ACC for lc thresholds for algebraically integrable foliations}, arXiv:2307.07157.


\bibitem[Dru21]{Dru21} S. Druel, \textit{Codimension 1 foliations with numerically trivial canonical class on singular spaces}. Duke Math. J., \textbf{170} (2021), no. 1, 95--203.



\bibitem[HL23]{HL23} C. D. Hacon and J. Liu, \textit{Existence of flips for generalized lc pairs}, Camb. J. Math. \textbf{11} (2023), no. 4, 795--828.  


\bibitem[HMX14]{HMX14} C. D. Hacon, J. M\textsuperscript{c}Kernan, and C. Xu, \textit{ACC for log canonical thresholds}, Ann. of Math. \textbf{180} (2014), no. 2, 523--571.


\bibitem[HLQ21]{HLQ21} J.~Han, Z.~Li, and L.~Qi, \textit{ACC for log canonical threshold polytopes}, Amer. J. Math. \textbf{143} (2021), no. 3, 681--714.


\bibitem[HL22]{HL22} J. Han and J. Liu, \textit{On termination of flips and exceptionally non-canonical singularities}, arXiv:2209.13122. To appear in Geom. Topol.

\bibitem[HLS19]{HLS19} J. Han, J. Liu, and V. V. Shokurov, \textit{ACC for minimal log discrepancies of exceptional singularities}, arXiv:1903.04338.

\bibitem[Kol21]{Kol21} J. Koll\'ar, \textit{Relative MMP without $\mathbb Q$-factoriality}, Electron. Res. Arch.  \textbf{29} (2021), no. 5, 3193--3203. 

\bibitem[Kol$^+$92]{Kol+92} J.~Koll\'{a}r \'{e}t al., \textit{Flip and abundance for algebraic threefolds}. Ast\'{e}risque no. \textbf{211} (1992).


\bibitem[KM98]{KM98} J. Koll\'{a}r and S. Mori, \textit{Birational geometry of algebraic varieties}, Cambridge Tracts in Math. \textbf{134} (1998), Cambridge Univ. Press.

\bibitem[Liu18]{Liu18} J. Liu, \textit{Toward the equivalence of the ACC for a-log canonical thresholds and the ACC for minimal log discrepancies}, arXiv:1809.04839.

\bibitem[LLM23]{LLM23} J. Liu, Y. Luo, and F. Meng, \textit{On global ACC for foliated threefolds}, Trans. Amer. Math. Soc. \textbf{376} (2023), no. 12, 8939--8972.

\bibitem[LMX24a]{LMX24a} J. Liu, F. Meng, and L. Xie, \textit{Complements, index theorem, and minimal log discrepancies of foliated surface singularities}, Eur. J. Math. \textbf{10} (2024), no. 6.

\bibitem[LMX24b]{LMX24b} J. Liu, F. Meng, and L. Xie \textit{Minimal model program for algebraically integrable foliations on klt varieties}, arXiv:2404.01559.

\bibitem[McQ08]{McQ08} M. McQuillan, \textit{Canonical models of foliations}, Pure Appl. Math. Q. \textbf{4} (2008), no. 3, Special Issue: In honor of Fedor Bogomolov, Part 2, 877--1012.

\bibitem[Miy87]{Miy87} Y. Miyaoka, \textit{Deformations of a morphism along a foliation and applications}, Algebraic geometry, Bowdoin, Proc. Sympos. Pure Math. \textbf{46} (1985) (Brunswick, Maine, 1985), Amer. Math. Soc., Providence, RI (1987), 245--268.

\bibitem[Nak16]{Nak16} Y.~Nakamura, \textit{On minimal log discrepancies on varieties with fixed Gorenstein index}. Michigan Math. J. \textbf{65} (2016), no. 1, 165--187.

\bibitem[Sho92]{Sho92} V.V.~Shokurov, \textit{Threefold log flips}, With an appendix in English by Y. Kawamata, Izv. Ross. Akad. Nauk Ser. Mat. \textbf{56} (1992), no. 1, 105--203 (Appendix by Y. Kawamata).

\bibitem[Spi20]{Spi20} C. Spicer, \textit{Higher dimensional foliated Mori theory}, Compos. Math. \textbf{156} (2020), no. 1, 1--38.

\bibitem[SS22]{SS22} C. Spicer and R. Svaldi, \textit{Local and global applications of the Minimal Model Program for co-rank 1 foliations on threefolds}, J. Eur. Math. Soc. \textbf{24} (2022), no. 11, 3969--4025.

\bibitem[SS23]{SS23} C. Spicer and R. Svaldi, \textit{Effective generation for foliated surfaces: Results and applications}, J. Reine Angew. Math. (2023), no. 795, 45--84.

\end{thebibliography}
\end{document}